\documentclass[10pt]{article}
\usepackage{latexsym}
\usepackage{amsfonts}
\usepackage{graphicx}
\usepackage{subfigure}
\usepackage{color}
\usepackage[active]{srcltx} %SRC Specials for DVI Searching
\usepackage{amssymb}
\usepackage[font=small,labelfont=bf]{caption}
\usepackage{hyperref}
\usepackage{soul}
\usepackage{multirow}

\newcommand{\ext}{\raise1pt\hbox{$\textstyle\bigwedge$}}

\title{Anisotropy in 2D Discrete Exterior Calculus}

\author
{Humberto Esqueda\footnote{Centro de Investigaci\'{o}n en Matem\'{a}ticas, Calle Jalisco s/n, 
Guanajuato, GTO 36023, M\'exico. Email: esqueda,rherrera,botello@cimat.mx
}, Rafael Herrera$^*$, Salvador Botello$^*$, \newline Carlos Valero\footnote{Departamento de 
Matem\'aticas, Universidad de Guanajuato, 
Guanajuato, GTO 36023, M\'exico. Email: valerocar@gmail.com}}

\date{\today} % delete this line to display the current date

\begin{document}

\maketitle

{

\abstract{

We present a local formulation for 2D Discrete Exterior Calculus (DEC) similar to that of the Finite Element Method (FEM), which allows 
a natural treatment of material heterogeneity (element by element). It also allows us to deduce, in a robust manner, anisotropic fluxes and
the DEC discretization of the pullback of 1-forms by the anisotropy tensor, i.e. we deduce how 
the anisotropy tensor acts on primal 1-forms. Due to the local formulation, the computational cost of DEC is similar to that of the Finite Element Method with Linear interpolations functions (FEML). The numerical 
DEC solutions to the anisotropic Poisson equation show numerical convergence, are very close to those of FEML on fine meshes 
and are slightly better than those of FEML on coarse meshes. 

}

}

\section{Introduction}

The theory of Discrete Exterior Calculus (DEC) is a relatively recent discretization  \cite{HiraniThesis} of the classical theory of 
Exterior Differential Calculus, a theory developed by 
E. Cartan \cite{Cartan} which has been a fundamental tool in Differential Geometry and Topology for over a century. 
The aim of DEC is to solve partial differential equations preserving their geometrical and physical
features as much as possible.
There are only a few papers about implementions of DEC to solve certain PDEs, such as 
the Darcy flow and Poisson's equation \cite{Hirani_K_N}, 
the Navier-Stokes equations \cite{Mohamedetal},
the simulation of  elasticity, plasticity and failure of isotropic materials \cite{Dassiosetal}, some comparisons with 
the finite differences and finite volume methods on regular flat meshes \cite{Griebel_R_S}, as well as applications in digital geometry processing \cite{Craneetal}.

In this paper, we describe a local formulation of DEC which is reminiscent of that of the Finite Element Method  (FEM) since,
once the local systems of equations have been established, they can be assembled into a global linear system.
This local formulation is also efficient and helpful in understanding various features of DEC that can otherwise remain unclear while dealing 
with an entire mesh. We will, therefore, take a local approach when recalling all the objects required by DEC \cite{Esqueda1}.
Our main results are the following: 
\begin{itemize}
\item We develop a local formulation of DEC analogous to that of FEM, which
allows a natural treatment of heterogeneous material properties assigned to subdomains (element by element) and 
eliminates the need of dealing with it through ad hoc modifications of the global discrete Hodge star operator. 

\item Guided by the local formulation,
we also deduce a natural way to approximate the flux/gradient-vector of a discretized function, as well as the anisotropic flux vector.
We carry out a comparison of the formulas defining the flux in both DEC and Finite Element Method with linear interpolation functions (FEML).

\item From the local formulation, we deduce the local DEC-discretization of the anisotropic Poisson equation. More precisely, in Exterior Differential Calculus the anisotropy tensor acts by {\em pullback} on the differential of the unknown function. Here, we deduce how the anisotropy tensor acts on primal 1-forms.
We also carry out an algebraic comparison of the DEC and FEML local formulations of the anisotropic Poisson equation.

\item We present three numerical examples of the approximate solutions to the stationary anisotropic Poisson equation on different domains using DEC and FEML. The numerical examples show numerical convergence and a competitive performance 
of DEC, as well as a computational cost similar to that of FEML. In fact, the numerical solutions with both methods on fine meshes are identical, and
DEC shows a slightly better performance than FEML on coarse meshes.
\end{itemize}

The paper is organized as follows.
In Section \ref{sec: preliminaries}, we describe the local versions of the discrete derivative operator,
the dual mesh and the discrete Hodge star operator. 
In Section \ref{sec: anisotropy}, we deduce the natural way of computing flux vectors in DEC (which turns out to be equivalent to the FEML procedure), as well as the anisotropic flux vectors. 
In Section \ref{sec: comparison}, 
we present the local DEC formulation of the 2D anisotropic Poisson equation
and compare it with the local system of FEML, proving that the diffusion terms are identical while the source terms are discretized differently due to a different area-weight assignment for the nodes. 
In Section \ref{sec: extension}, we re-examine some of the local DEC quantities.
In Section \ref{sec: examples}, we present and compare numerical examples of DEC and FEML approximate solutions to the 2D anisotropic 
Poisson equation on different domains with meshes of various resolutions. 
In Section \ref{sec: conclusions} we summarize the contributions of this paper.

{\em Acknowledgements}. The second named author was partially supported by a CONACyT grant,
and would like to thank 
the International Centre for Numerical Methods in Engineering (CIMNE) and the University of Swansea
for their hospitality. We gratefully acknowledge the support of
NVIDIA Corporation with the donation of the Titan X Pascal GPU used for this research.

\section{Preliminaries on DEC from a local viewpoint}\label{sec: preliminaries}

Let us consider a primal mesh made up of a single (positively oriented) triangle.
    \begin{center}
    \includegraphics[width=0.4\textwidth]{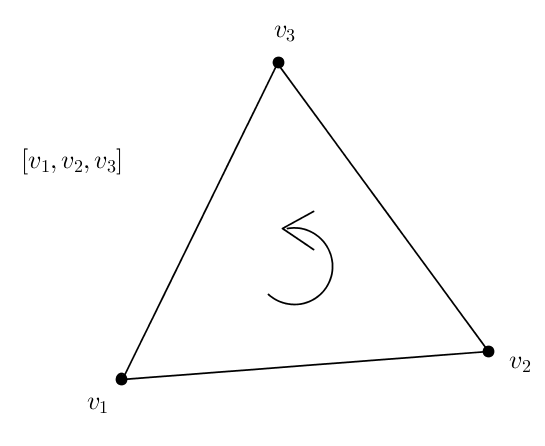}%
    \captionof{figure}{Triangle $[v_1,v_2,v_3]$.}
    \end{center}
\subsection{Boundary operator}
There is a well known boundary operator 
\begin{equation}
\partial [v_1,v_2,v_3]=[v_2,v_3]-[v_1,v_3]+[v_1,v_2],\label{eq: boundary 2-1}
\end{equation}
which describes the boundary of the triangle as an alternated sum of its oriented edges $[v_1,v_2]$, $[v_2,v_3]$ and $[v_3,v_1]$. 
Similarly, one can compute the 
boundary of each edge
\begin{eqnarray}
\partial [v_1,v_2] &=& [v_2]-[v_1], \nonumber\\ 
\partial [v_2,v_3] &=& [v_3]-[v_2], \label{eq: boundary 1-0}\\ 
\partial [v_3,v_1] &=& [v_1]-[v_3]. \nonumber
\end{eqnarray}
If we consider 
\begin{itemize} 
\item the symbol $[v_1,v_2,v_3]$ as a basis vector of a 1-dimensional vector space,
\item the symbols $[v_1,v_2]$, $[v_2,v_3]$, $[v_3,v_1]$ as an ordered basis of a 3-dimensional vector space,  
\item the symbols $[v_1]$, $[v_2]$, $[v_3]$ as an ordered basis of a 3-dimensional vector space,
\end{itemize}
then the map (\ref{eq: boundary 2-1}), which sends the oriented triangle to a sum of its oriented edges, is represented by the matrix
\[\left(\begin{array}{r}
1 \\ 
1 \\ 
1
\end{array}\right) ,\]
while the map (\ref{eq: boundary 1-0}), which sends the oriented edges to sums of their oriented vertices, is represented by the matrix
\[\left(\begin{array}{rrr}
-1 & 0 & 1 \\ 
1 & -1 & 0 \\ 
0 & 1 & -1
\end{array} \right).\]

\subsection{Discrete derivative}\label{subsec: discrete derivative}
It has been argued that the DEC discretization of the differential of a function is given by the transpose of the matrix of the 
boundary operator on edges (see \cite{HiraniThesis, Esqueda1}). More precisely, suppose we have a function discretized by its values at the vertices
\[f\sim \left(\begin{array}{l}
f_1\\
f_2\\
f_3
              \end{array}
\right).\]
Its discrete derivative, according to DEC, is 
\begin{eqnarray*}
\left(\begin{array}{rrr}
-1 & 0 & 1 \\ 
1 & -1 & 0 \\ 
0 & 1 & -1
\end{array} \right)^T
\left(\begin{array}{l}
f_1\\
f_2\\
f_3
              \end{array}
\right)
&=&
\left(\begin{array}{rrr}
-1 & 1 & 0 \\ 
0 & -1 & 1 \\ 
1 & 0 & -1
\end{array} \right)
\left(\begin{array}{l}
f_1\\
f_2\\
f_3
              \end{array}
\right)\\
&=&
\left(\begin{array}{l}
f_2-f_1\\
f_3-f_2\\
f_1-f_3
              \end{array}
\right) .
\end{eqnarray*}
Indeed, such differences are rough approximations of the directional derivatives of $f$. For instance, $f_2-f_1$ is a rough approximation of the directional derivative of $f$ at $v_1$ in the direction of the vector $v_2-v_1$, i.e.
\[f_2-f_1 \approx df_{v_1}(v_2-v_1) .\]
It is precisely in this sense that, according to DEC, 
\begin{itemize}
\item the value $f_2-f_1$ is assigned to the edge $[v_1,v_2]$, 
\item the value $f_3-f_2$ is assigned to the edge $[v_2,v_3]$, 
\item and the value $f_1-f_3$ is assigned to the edge $[v_3,v_1]$.
\end{itemize}
Let
\[D_0:=\left(\begin{array}{rrr}
-1 & 1 & 0 \\ 
0 & -1 & 1 \\ 
1 & 0 & -1
\end{array} \right).
\]

\subsection{Dual mesh}

The dual mesh of the primal mesh consisting of a single triangle is constructured as follows:
\begin{itemize}
 \item To the 2-dimensional triangular face $[v_1,v_2,v_3]$ will correspond the 0-dimensional point given by the circumcenter $c$ of the triangle.
     \begin{center}
    \includegraphics[width=0.4\textwidth]{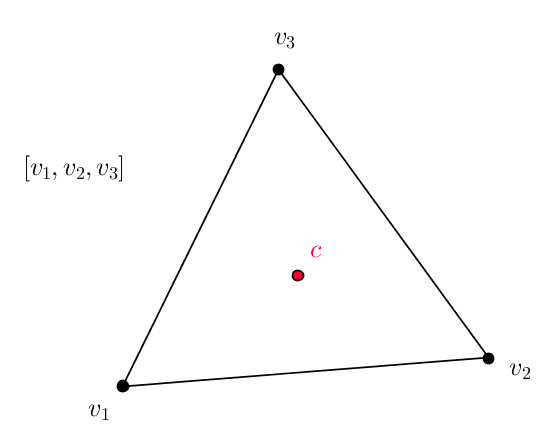}%
    \captionof{figure}{Circumcenter $c$ of the triangle $[v_1,v_2,v_3]$.}
    \end{center}

 \item To the 1-dimensional edge $[v_1,v_2]$ will correspond the 1-dimensional straight line segment $[p_1,c]$ joining the midpoint $p_1$ of the edge $[v_1,v_2]$ to the circumcenter $c$. 
 Similarly for the other edges.
     \begin{center}
    \includegraphics[width=0.4\textwidth]{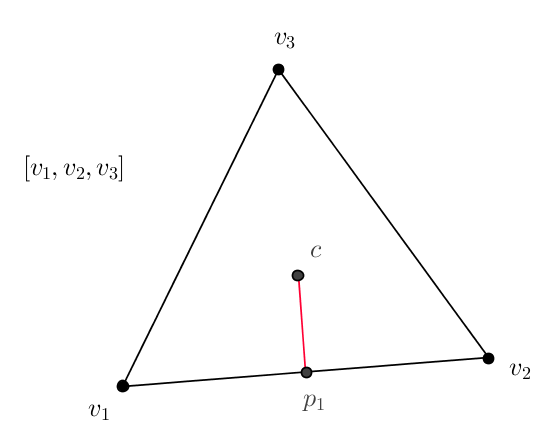}%
    \captionof{figure}{Dual segment $[p_1,c]$ of the edge $[v_1,v_2]$.}
    \end{center}

 \item To the 0-dimensional vertex/node $[v_1]$ will correspond the 2-dimensional quadrilateral $[v_1,p_1,c,p_3]$.
    \begin{center}
    \includegraphics[width=0.4\textwidth]{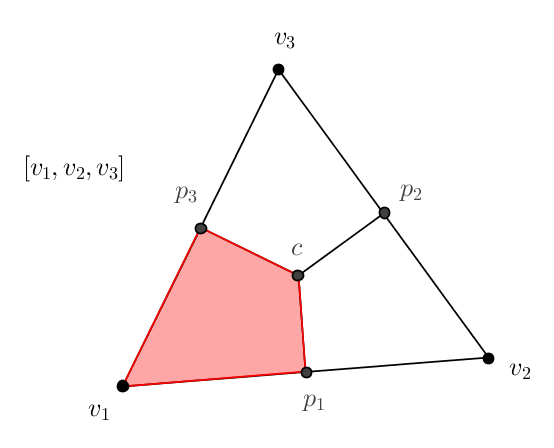}%
    \captionof{figure}{Dual quadrilateral  $[v_1,p_1,c,p_3]$ of the vertex $[v_1]$.}
    \end{center}

\end{itemize}

\subsection{Discrete Hodge star}

For the Poisson equation in 2D, we need two matrices: one relating original edges to dual edges, and another relating vertices to dual cells.
\begin{itemize}
\item The discrete Hodge star map $M_1$ applied to the discrete differential of a discretized function $f\sim (f_1,f_2,f_3)$ 
is given as follows:
\begin{itemize}
 \item the value $f_2-f_1$ assigned to the edge $[v_1,v_2]$ is changed to the new value 
 \[{{\rm length}[p_1,c]\over {\rm length}[v_1,v_2]}(f_2-f_1)\]
 assigned to the segment $[p_1,c]$;
 \item the value $f_3-f_2$ assigned to the edge $[v_2,v_3]$ is changed to the new value 
 \[{{\rm length}[p_2,c]\over {\rm length}[v_2,v_3]}(f_3-f_2)\]
 assigned to the segment $[p_2,c]$;
 \item the value $f_1-f_3$ assigned to the edge $[v_3,v_1]$ is changed to the new value 
 \[{{\rm length}[p_3,c]\over {\rm length}[v_3,v_1]}(f_1-f_3)\]
 assigned to the edge $[p_3,c]$.
\end{itemize}
In other words, 
\[M_1=\left(\begin{array}{ccc}
{{\rm length}[p_1,c]\over {\rm length}[v_1,v_2]} & 0 & 0 \\ 
0 & {{\rm length}[p_2,c]\over {\rm length}[v_2,v_3]} & 0 \\ 
0 & 0 & {{\rm length}[p_3,c]\over {\rm length}[v_3,v_1]}
\end{array} \right).\]

\item Similarly, the discrete Hodge star map $M_0$ on values on vertices is given as follows
\begin{itemize}
 \item the value $f_1$ assigned to the vertex $[v_1]$ is changed to the new value 
 \[{\rm Area}[v_1,p_1,c,p_3]f_1\]
 assigned to the quadrilateral $[v_1,p_1,c,p_3]$;
 \item the value $f_2$ assigned to the vertex $[v_2]$ is changed to the new value
 \[{\rm Area}[v_2,p_2,c,p_1]f_2\]
 assigned to the quadrilateral $[v_2,p_2,c,p_1]$;
 \item the value $f_3$ assigned to the vertex $[v_3]$ is changed to the new value
 \[{\rm Area}[v_3,p_3,c,p_2]f_2\]
 assigned to the quadrilateral $[v_3,p_3,c,p_2]$.
\end{itemize}
In other words,
\[M_0=\left(\begin{array}{ccc}
{\rm Area}[v_1,p_1,c,p_3] & 0 & 0 \\ 
0 & {\rm Area}[v_2,p_2,c,p_1] & 0 \\ 
0 & 0 & {\rm Area}[v_3,p_3,c,p_2]
\end{array} \right).\]
\end{itemize}

\section{Flux and anisotropy}\label{sec: anisotropy}
In this section, we deduce the DEC formulae for the local flux, the local anisotropic flux and the local anisotropy operator 
for primal 1-forms.

\subsection{The flux in local DEC}\label{subsec: local DEC flux}
We wish to find a natural construction for the discrete flux (discrete gradient vector) of a discrete function.
Recall from Vector Calculus that the directional derivative of a differentiable function $f:\mathbb{R}^2\longrightarrow \mathbb{R}$ at a point $p\in \mathbb{R}^2$ in the direction of $w\in\mathbb{R}^2$ is defined by
\begin{eqnarray*}
df_p(w)&:=&\lim_{t\rightarrow 0}{f(p+tw)-f(p)\over t}\\
&=&\nabla f(p)\cdot w.
\end{eqnarray*}
Thus, we have three Vector Calculus identities
\begin{eqnarray*}
 df_{v_1}(v_2-v_1) &=& \nabla f(v_1)\cdot (v_2-v_1)  ,\\
 df_{v_2}(v_3-v_2) &=& \nabla f(v_2)\cdot (v_3-v_2),\\
 df_{v_3}(v_1-v_3) &=& \nabla f(v_3)\cdot (v_1-v_3) .
\end{eqnarray*}
As in subsection \ref{subsec: discrete derivative}, the rough approximations to directional derivatives of a function $f$ in the directions of the (oriented) edges are given as follows
\begin{eqnarray*}
 df_{v_1}(v_2-v_1) &\approx& f_2-f_1  ,\\
 df_{v_2}(v_3-v_2) &\approx& f_3-f_2  ,\\
 df_{v_3}(v_1-v_3) &\approx& f_1-f_3  .
\end{eqnarray*}
Thus, if we want to find a discrete gradient vector $W_1$ of $f$
at the point $v_1$,  we need to solve the equations of approximations
\begin{eqnarray}
W_1\cdot (v_2-v_1) &=& f_2-f_1 \label{eq: W1-1}\\
W_1\cdot (v_3-v_1) &=& f_3-f_1.\label{eq: W1-2}
\end{eqnarray} 
If
\begin{eqnarray*}
v_1&=&(x_1,y_1),\\
v_2&=&(x_2,y_2),\\
v_3&=&(x_3,y_3),
\end{eqnarray*}
then
\begin{eqnarray*}
W_1
&=&
\left({f_1y_2-f_1y_3-f_2y_1+f_2y_3+f_3y_1-f_3y_2\over x_1y_2-x_1y_3-x_2y_1+x_2y_3+x_3y_1-x_3y_2}
,  -{f_1x_2-f_1x_3-f_2x_1+f_2x_3+f_3x_1-f_3x_2\over x_1y_2-x_1y_3-x_2y_1+x_2y_3+x_3y_1-x_3y_2}\right)^T
\end{eqnarray*}

Now, if we were to find 
a discrete gradient vector $W_2$ of $f$
at the point $v_2$,  we need to solve the equations
% the approximate $\nabla f(v_2)$ we need to find $W_2$ such that
\begin{eqnarray}
W_2\cdot (v_1-v_2) &=& f_1-f_2\label{eq: W2-1}\\
W_2\cdot (v_3-v_2) &=& f_3-f_2.\nonumber %\label{eq: W2-2}
\end{eqnarray} 
The vectors $W_2$ solving these equations is actually equal to $W_1$. 
Indeed, consider
\begin{eqnarray*}
f_3-f_1 &=& W_1\cdot (v_3-v_1) \\
&=& W_1\cdot (v_3-v_2+v_2-v_1)\\
&=& W_1\cdot (v_3-v_2)+W_1\cdot (v_2-v_1)\\
&=& W_1\cdot (v_3-v_2)+f_2-f_1,
\end{eqnarray*} 
so that
\begin{equation}
 W_1\cdot (v_3-v_2) = f_3-f_2  .\label{eq: W1-3}
\end{equation} 
Thus, adding up (\ref{eq: W1-1}) and (\ref{eq: W2-1}) we get
\begin{equation}
(W_1-W_2)\cdot (v_2-v_1)=0.\label{eq: W1-W2-1}
\end{equation}
Subtracting  (\ref{eq: W1-2}) from (\ref{eq: W1-3}) we get
\begin{equation}
(W_1-W_2)\cdot (v_3-v_2)=0.\label{eq: W1-W2-2}
\end{equation}
Since $v_2-v_1$ and $v_3-v_2$ are linearly independent 
and the two inner products in (\ref{eq: W1-W2-1}) and (\ref{eq: W1-W2-2}) vanish,
\[W_1-W_2=0.\]

Analogously, the corresponding gradient vector $W_3$ of $f$ at the vertex $v_3$ is equal to $W_1$.
This means that the three approximate gradient vectors at the three vertices coincide. 
Let us call this unique vector $W$.
Note that discrete flux $W$ satisfies 
\begin{eqnarray*}
W\cdot (v_2-v_1) &=& f_2-f_1,\\
W\cdot (v_3-v_1) &=& f_3-f_1,\\
W\cdot (v_3-v_2) &=& f_3-f_2  .
\end{eqnarray*}
This means that the primal 1-form discretizing $df$ can be obtained by the dot products of the discrete flux $W$ with the 
vectors of the triangle's edges.

\vspace{.1in}

{\bf Remark}.
More generally, we can see that {\em any vector which is constant on the triangle, naturally gives a primal 1-form on the edges of the triangle 
by means of its dot products with the triangle's edge-vectors}.

\subsubsection{Comparison of DEC and FEML local fluxes}

The local flux (gradient) of $f$ in FEML is given by
\[
\left(
\begin{array}{ccc}
{\partial N_1\over\partial x} & {\partial N_2\over\partial x} & {\partial N_3\over\partial x} \\ 
{\partial N_1\over\partial y} & {\partial N_2\over\partial y} & {\partial N_3\over\partial y}
\end{array} 
\right)
\left(
\begin{array}{c}
f_1 \\ 
f_2 \\ 
f_3
\end{array} 
\right),
\]
where
\begin{eqnarray*}
N_1 &=& {1\over 2A}[(y_2-y_3)x+(x_3-x_2)y+x_2y_3-x_3y_2] ,\\
N_2 &=& {1\over 2A}[(y_3-y_1)x+(x_1-x_3)y+x_3y_1-x_1y_3] ,\\
N_3 &=& {1\over 2A}[(y_1-y_2)x+(x_2-x_1)y+x_1y_2-x_2y_1] ,
\end{eqnarray*}
and
\begin{eqnarray*}
A &=& {1\over 2} [(x_2y_3-x_3y_2) -(x_1y_3-x_3y_1)+(x_1y_2-x_2y_1)]
\end{eqnarray*}
is the area of the triangle.
Explicitly
\begin{eqnarray*}
\left(
\begin{array}{ccc}
{\partial N_1\over\partial x} & {\partial N_2\over\partial x} & {\partial N_3\over\partial x} \\ 
{\partial N_1\over\partial y} & {\partial N_2\over\partial y} & {\partial N_3\over\partial y}
\end{array} 
\right)
&=&
{1\over 2A}
\left(
\begin{array}{ccc}
y_2-y_3 & y_3-y_1 & y_1-y_2 \\ 
x_3-x_2 & x_1-x_3 & x_2-x_1
\end{array} 
\right)  
\end{eqnarray*}
so that the FEML flux is given by
\[
\left(
{[(y_2-y_3)f_1+(y_3-y_1)f_2+(y_1-y_2)f_3]\over 2A} ,
{[(x_3-x_2)f_1+(x_1-x_3)f_2+(x_2-x_1)f_3]\over 2A}
 \right)^T,
\]
and we can see that its formula coincides with that of the DEC flux.

\subsection{The anisotropic flux vector in local DEC} \label{subsec: discrete anisotropic flux}

We will now discuss how to discretize anisotropy in 2D DEC.
Let $K$ denote the symmetric anisotropy tensor 
\[K=\left(\begin{array}{cc}
k_{11} & k_{12} \\ 
k_{12} & k_{22}
\end{array} \right)\]
and recall the anisotropic Poisson equation
\[-\nabla\cdot (K\, \nabla f) = q.\]
As in Subsection \ref{subsec: local DEC flux}, we wish to find a vector $W'$ which will play the role of a discrete version of the {\em anisotropic flux vector} $K \nabla f$.

First observe that, since $K$ is symmetric, for any $w\in\mathbb{R}^2$
\begin{eqnarray*}
(K\nabla f(p))\cdot w 
   &=& 
  \nabla f(p)\cdot (K^Tw)\\
   &=& 
  \nabla f(p)\cdot (Kw)\\
   &=&
  df_p(Kw)\\
   &=&
  (df_p\circ K)(w)\\
   &=:&
  (K^*df_p)(w),
\end{eqnarray*}
where $K^*df_p$ is called the {\em pullback} of $df_p$ by $K$. 
These identities mean that in order to discretize the anisotropic flux we need to understand the discretization
of the linear functional $df_p\circ K$. 
Let us suppose that {\em $K$ is constant on our triangle}.
As before, we have three natural vectors on the triangle, 
\begin{eqnarray*}
w_1 &=& v_2-v_1,\\
w_2 &=& v_3-v_2,\\
w_3 &=& v_1-v_3.
\end{eqnarray*}
Given the vector $Kw_1$, we have the option to write it down as a linear combination of two of the three aforementioned vectors.
Since $w_1$ is being used already, we use the other two vectors, i.e.
\begin{eqnarray*}
Kw_1 &=& \lambda_1 w_2 + \mu_1 w_3,
\end{eqnarray*}
for some $\lambda_1,\mu_1\in\mathbb{R}$.
Similarly,
\begin{eqnarray*}
Kw_2 &=& \lambda_2 w_3 + \mu_2 w_1,\\
Kw_3 &=& \lambda_3 w_1 + \mu_3 w_2,
\end{eqnarray*}
for some $\lambda_2,\mu_2,\lambda_3,\mu_3\in\mathbb{R}$.
These equations can be solved for $\lambda_1,\lambda_2,\lambda_3,\mu_1,\mu_2,\mu_3$.
Now
\begin{eqnarray*}
\left(K\nabla f(v_3)\right)\cdot w_1 
   &=&
  df_{v_3}(Kw_1)\\
   &=&
  df_{v_3}(\lambda_1 w_2 + \mu_1 w_3) \\
   &=&
  \lambda_1 df_{v_3}(w_2) + \mu_1 df_{v_1}(w_3) \\
   &=&
  \lambda_1 df_{v_3}(v_3-v_2) + \mu_1 df_{v_3}(v_1-v_3) \\
   &=&
  \lambda_1 df_{v_3}(-(v_2-v_3)) + \mu_1 df_{v_3}(v_1-v_3) \\
   &=&
  -\lambda_1 df_{v_3}(v_2-v_3) + \mu_1 df_{v_3}(v_1-v_3) .
\end{eqnarray*}
Similarly, for the vectors $w_2$ and $w_3$ we have the identities 
\begin{eqnarray*}
\left(K\nabla f(v_1)\right)\cdot w_2 &=& \lambda_2 df_{v_1}(w_3) + \mu_2 df_{v_1}(w_1),\\
\left(K\nabla f(v_2)\right)\cdot w_3  &=& \lambda_3 df_{v_3}(w_1) +\mu_3 df_{v_3}(w_2).
\end{eqnarray*}
These equations lead to the three equations of approximations
\begin{eqnarray}
W'\cdot w_1
   & = &
  \lambda_1 (f_3-f_2) + \mu_1 (f_1-f_3), \nonumber \\
W'\cdot w_2
   & = &
  \lambda_2 (f_1-f_3) + \mu_2 (f_2-f_1), \label{eq: system equations 1 for anisotropic flux vector}\\
W'\cdot w_3
   & = &
  \lambda_3 (f_2-f_1) + \mu_3 (f_3-f_2). \nonumber
\end{eqnarray}
where $W'$ is the vector that should approximate $K\nabla f(v_3)$, $K\nabla f(v_3)$ and $K\nabla f(v_3)$.
Thus, in order to find the discrete version $W'$ of the anisotropic flux vector $K\nabla f$ on the triangle, we need to solve 
the system (\ref{eq: system equations 1 for anisotropic flux vector})

The system (\ref{eq: system equations 1 for anisotropic flux vector}) has a unique solution. Indeed, since
\[w_1+w_2+w_3=0,\] 
then
\[Kw_1+Kw_2+Kw_3=0,\] 
i.e.
\begin{eqnarray*}
(\lambda_3+\mu_2-\lambda_2-\mu_1) w_1 + (\lambda_1+\mu_3-\lambda_2-\mu_1) w_2
    =0.
\end{eqnarray*}
Since $w_1$ and $w_2$ are linearly independent
\begin{eqnarray*}
\lambda_3+\mu_2-\lambda_2-\mu_1&=&0\\
\lambda_1+\mu_3-\lambda_2-\mu_1&=&0.
\end{eqnarray*}
i.e.
\begin{eqnarray*}
\lambda_3&=&\lambda_2+\mu_1-\mu_2\\
\mu_3&=&-\lambda_1+\lambda_2+\mu_1.
\end{eqnarray*}
Thus, making the appropriate substitutions,  we see that
the third equation in (\ref{eq: system equations 1 for anisotropic flux vector}) is dependent on the first two independent equations,
and there is a unique vector $W'$ that solves the system.

For the sake of completeness, the values of the parameters are:
\begin{eqnarray*}
\lambda_1
   &=&
  {[k_{11}(x_2-x_1)+k_{12}(y_2-y_1)](y_1-y_3) -[k_{12}(x_2-x_1)+k_{22}(y_2-y_1)](x_1-x_3)
  \over (x_3-x_2)(y_1-y_3)-(x_1-x_3)(y_3-y_2)}\\ 
   &=&
   -{J(w_3)\cdot K(w_1) \over 2A},
   \\
\mu_1
   &=&
  -{[[k_{11}(x_2-x_1)+k_{12}(y_2-y_1)](y_3-y_2)-[k_{12}(x_2-x_1)+k_{22}(y_2-y_1)](x_3-x_2) 
  \over (x_3-x_2)(y_1-y_3)-(x_1-x_3)(y_3-y_2)}\\ 
   &=&
   {J(w_2)\cdot K(w_1) \over 2A},
   \\
\lambda_2
   &=&
  {[k_{11}(x_3-x_2)+k_{12}(y_3-y_2)](y_2-y_1) -[k_{12}(x_3-x_2)+k_{22}(y_3-y_2)](x_2-x_1)
  \over (x_1-x_3)(y_2-y_1)-(x_2-x_1)(y_1-y_3)}\\ 
   &=&
   -{J(w_1)\cdot K(w_2) \over 2A},
   \\
\mu_2
   &=&
  -{[[k_{11}(x_3-x_2)+k_{12}(y_3-y_2)](y_1-y_3)-[k_{12}(x_3-x_2)+k_{22}(y_3-y_2)](x_1-x_3) 
  \over (x_1-x_3)(y_2-y_1)-(x_2-x_1)(y_1-y_3)}\\ 
   &=&
   {J(w_3)\cdot K(w_2) \over 2A},
   \\
\lambda_3
   &=&
  {[k_{11}(x_1-x_3)+k_{12}(y_1-y_3)](y_3-y_2) -[k_{12}(x_1-x_3)+k_{22}(y_1-y_3)](x_3-x_2)
  \over (x_2-x_1)(y_3-y_2)-(x_3-x_2)(y_2-y_1)}\\ 
   &=&
   -{J(w_2)\cdot K(w_3) \over 2A},
   \\
\mu_3
   &=&
  -{[[k_{11}(x_1-x_3)+k_{12}(y_1-y_3)](y_2-y_1)-[k_{12}(x_1-x_3)+k_{22}(y_1-y_3)](x_2-x_1) 
  \over (x_2-x_1)(y_3-y_2)-(x_3-x_2)(y_2-y_1)}\\ 
   &=&
   {J(w_1)\cdot K(w_3) \over 2A},   
\end{eqnarray*}
where 
\begin{eqnarray*}
J
   &=&
  \left(\begin{array}{cc}
  0 & -1 \\ 
  1 & 0
  \end{array} \right)
\end{eqnarray*}
is the $90^\circ$ counter-clockwise rotation,
and
\begin{eqnarray*}
W'
   &=&
\left(
\begin{array}{c}
{k_{11}[f_1(y_2-y_3)+f_2(y_3-y_1)+f_3(y_1-y_2)]
+k_{12}[f_1(x_3-x_2)+f_2(x_1-x_3)+f_3(x_2-x_1)]
\over x_1y_2-x_1y_3-x_2y_1+x_2y_3+x_3y_1-x_3y_2} \\
{k_{12}[f_1(y_2-y_3)+f_2(y_3-y_1)+f_3(y_1-y_2)]
+k_{22}[f_1(x_3-x_2)+f_2(x_1-x_3)+f_3(x_2-x_1)]
\over x_1y_2-x_1y_3-x_2y_1+x_2y_3+x_3y_1-x_3y_2}
\end{array}
\right)\\
&=&
\left(\begin{array}{cc}
k_{11} & k_{12} \\ 
k_{12} & k_{22}
\end{array} \right)
\left(
\begin{array}{c}
{[(y_2-y_3)f_1+(y_3-y_1)f_2+(y_1-y_2)f_3]\over 2A} \\ 
{[(x_3-x_2)f_1+(x_1-x_3)f_2+(x_2-x_1)f_3]\over 2A}
\end{array} 
\right)
\end{eqnarray*}
which is, in fact, {\em the image under $K$ of the discrete isotropic flux} and shows the consistency of our {\em local} reasoning.  
Also observe that this formula is the same as that of the FEML anisotropic flux.

\subsection{Anisotropy on primal 1-forms}
The system (\ref{eq: system equations 1 for anisotropic flux vector}) can be rewritten in matrix form as follows
\begin{eqnarray}
\left(\begin{array}{c}
W'\cdot w_1 \\ 
W'\cdot w_2 \\ 
W'\cdot w_3
\end{array} \right)
 &=&
\left(
\begin{array}{ccc}
0 & \lambda_1 & \mu_1 \\ 
\mu_2 & 0 & \lambda_2 \\ 
\lambda_3 & \mu_3 & 0
\end{array} 
\right)
 \left(\begin{array}{l}
f_2-f_1\\
f_3-f_2\\
f_1-f_3
              \end{array}
\right) \nonumber\\
 &=&
\left(
\begin{array}{ccc}
0 & \lambda_1 & \mu_1 \\ 
\mu_2 & 0 & \lambda_2 \\ 
\lambda_3 & \mu_3 & 0
\end{array} 
\right)
\left(\begin{array}{rrr}
-1 & 1 & 0 \\ 
0 & -1 & 1 \\ 
1 & 0 & -1
\end{array} \right)
 \left(\begin{array}{l}
f_1\\
f_2\\
f_3
              \end{array}
\right) \nonumber\\
 &=&
\left(
\begin{array}{ccc}
0 & \lambda_1 & \mu_1 \\ 
\mu_2 & 0 & \lambda_2 \\ 
\lambda_3 & \mu_3 & 0
\end{array} 
\right)
D_0[f] .\label{eq: dual of W}
\end{eqnarray}
Recalling the Remark at the end of Subsection \ref{subsec: local DEC flux}, the matrix identity (\ref{eq: dual of W})
states that the primal 1-form dual to the anisotropic flux vector $W'$ is given by the local DEC discretization 
of the anisotropy tensor 
\begin{eqnarray*}
K^{DEC}
   &=&
\left(
\begin{array}{ccc}
0 & \lambda_1 & \mu_1 \\ 
\mu_2 & 0 & \lambda_2 \\ 
\lambda_3 & \mu_3 & 0
\end{array} 
\right)\\
   &=&
{1 \over 2A}
\left(
\begin{array}{ccc}
0 &    -J(w_3)\cdot K(w_1) & J(w_2)\cdot K(w_1) \\ 
J(w_3)\cdot K(w_2) & 0 & -J(w_1)\cdot K(w_2) \\ 
-J(w_2)\cdot K(w_3) & J(w_1)\cdot K(w_3) & 0
\end{array} 
\right).
\end{eqnarray*}
acting on the primal 1-form $D_0[f]$.

{\bf Remark}. The matrix $K^{DEC}$ is the local DEC discretization on primal 1-forms of the pullback operator 
$K^*$ on 1-forms. In this case, the discretization of
$K^*df:=df\circ K$.

\subsubsection{Geometric interpretation of the entries of $K^{DEC}$}
Let us examine $\lambda_1$ in the anisotropic case.
Consider the following figure
    \begin{center}
    \includegraphics[width=0.6\textwidth]{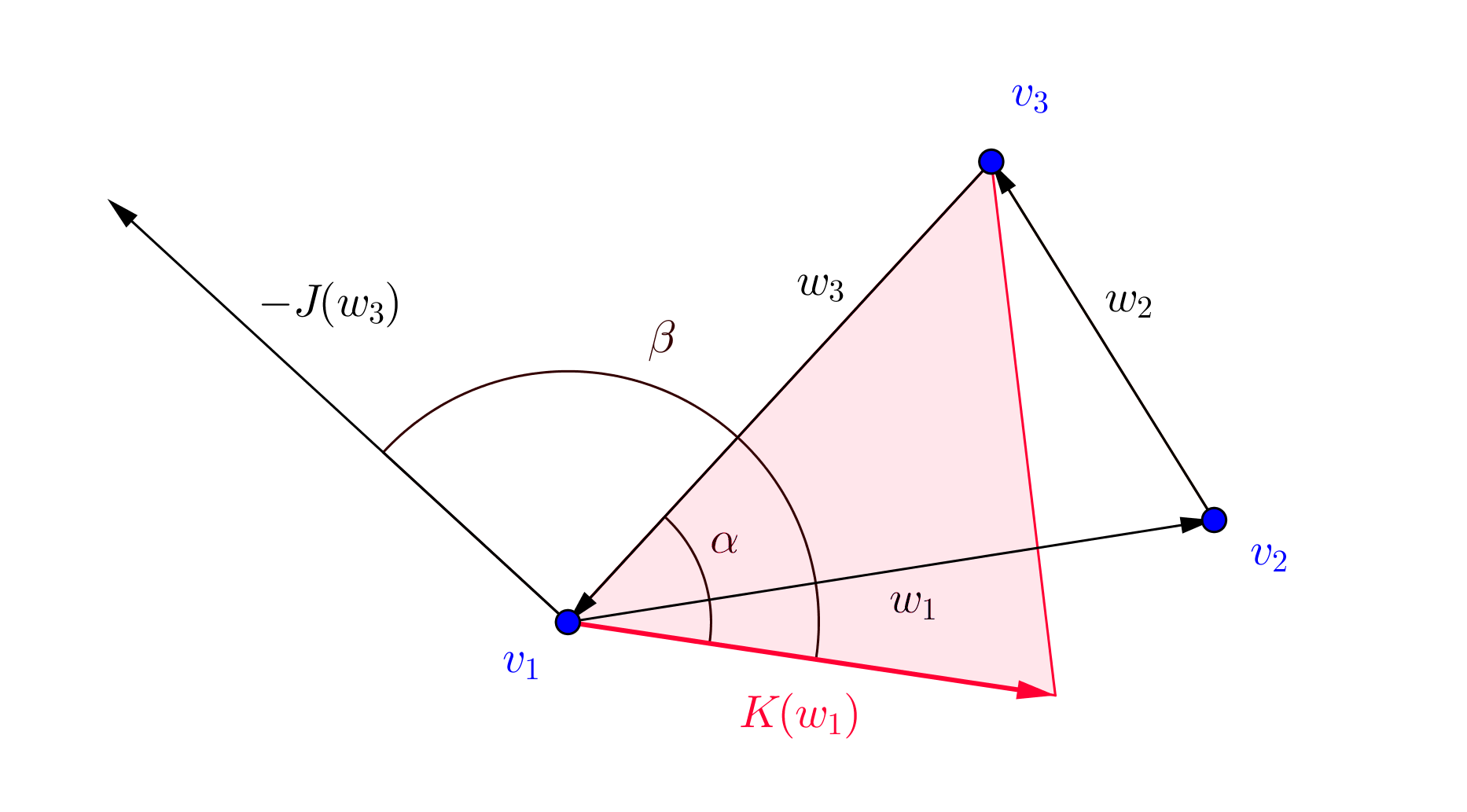}%
    \captionof{figure}{Geometric interpretation of the entries of the anisotropy tensor discretization $K^{ DEC}$}\label{fig: anisotropy triangle 3}
    \end{center}
We have
\begin{eqnarray*}
\lambda_1
   &= &   
 -{J(w_3)\cdot K(w_1) \over 2A}\\
   &= &   
 {1\over 2A} |-J(w_3)| |K(w_1)| \cos(\beta) \\
   &= &   
 {1\over 2A} |w_3| |K(w_1)| \cos(\alpha+\pi/2) \\
   &= &   
 -{1\over 2A} |w_3| |K(w_1)| \sin(\alpha) \\
   &=&
 -{A'\over A},
\end{eqnarray*}
where $A'$ is the area of the red triangle and
we have used a well known formula for the area of a triangle in terms of an inner angle. Thus
$\lambda_1$ is the negative of the quotient of the area $A'$ of the red triangle and the area $A$ of the original triangle.
The calculations for the other entries are similar.

\subsubsection{Isotropic case}
Now, let us assume $K=k\,{\rm Id}_{2\times 2}$ on the triangle.
The previous calculations show that
\[K^{DEC}=
k\left(\begin{array}{ccc}
0 & -1 & -1 \\ 
-1 & 0 & -1 \\ 
-1 & -1 & 0
\end{array} \right).\]
Note that, in this case,
\[K^{DEC} D_0= k\, D_0.\]

\section{2D anisotropic Poisson equation}\label{sec: comparison}
In this section, we describe the local DEC discretization of the 2D anisotropic Poisson equation and compare it to that of FEML.

\subsection{Local DEC discretization of the 2D anisotropic Poisson equation}
The anisotropic Poisson equation reads as follows
\[-\nabla\cdot (K\, \nabla f) = q,\]
where $f$ and $q$ are two functions on a certain domain in $\mathbb{R}^2$.
In terms of the exterior derivative $d$ and the Hodge star operator $\star$ it reads as follows
\[-\star d \star (K^*df) = q\]
where $K^*df:=df\circ K$
and $K=K^T$.
Following the discretization of the discretized divergence operator \cite{Esqueda1},
the corresponding local DEC discretization of the anisotropic Poisson equation is
\[-M_0^{-1} \, \left(-D_0^T\right)  \, M_1\, K^{DEC} \, D_0 \, [f] = [q], \]
or equivalently 
\begin{equation}
D_0^T  \, M_1\, K^{DEC}\, D_0 \, [f] = M_0 \, [q]. \label{eq: local DEC discretized anisotropic Poisson equation}
\end{equation} 
In order to simplify the notation,
consider the lengths and areas defined in the Figure \ref{fig: triangulo}.
    \begin{center}
    \includegraphics[width=0.35\textwidth]{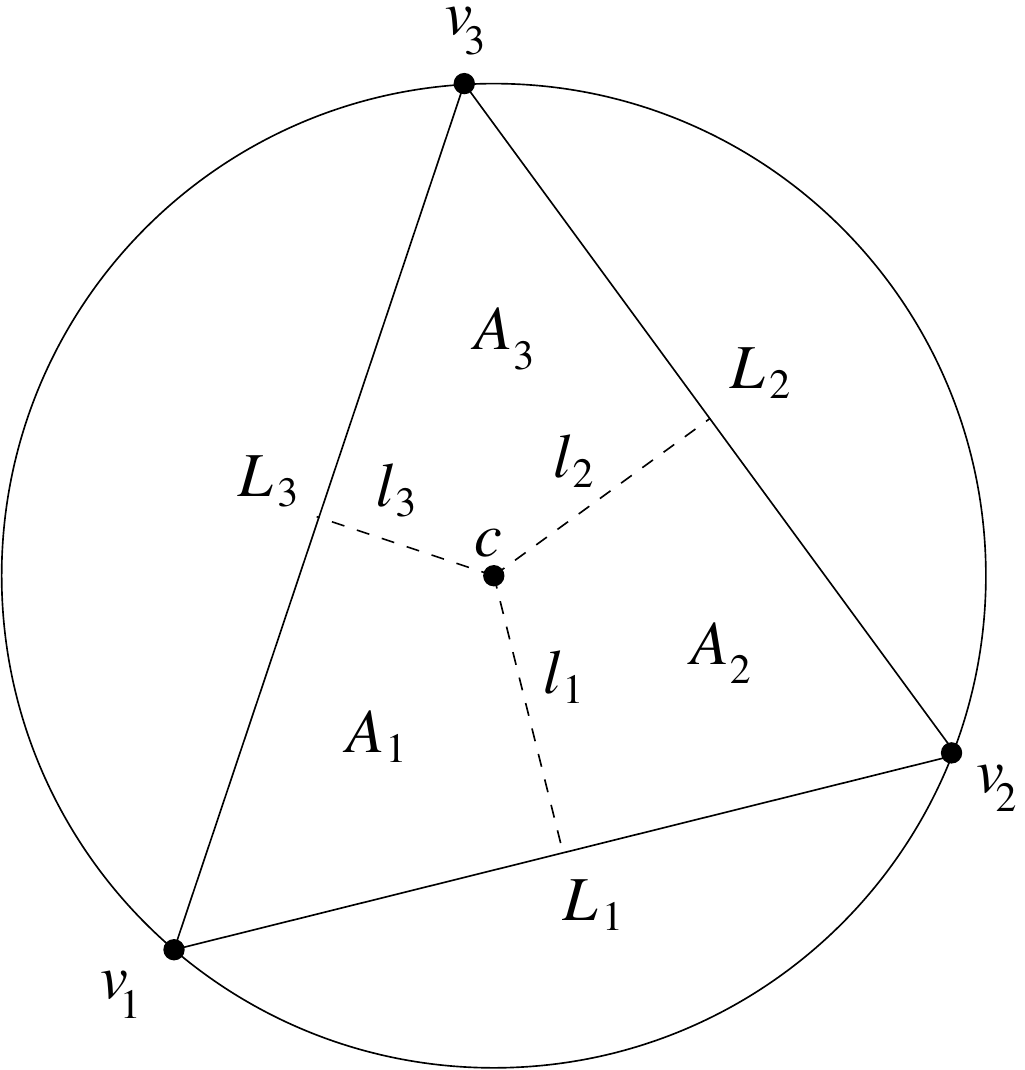}%
    \captionof{figure}{Triangle }\label{fig: triangulo}
    \end{center}
Now, the discretized equation (\ref{eq: local DEC discretized anisotropic Poisson equation}) looks as follows:
\[
\left(\begin{array}{ccc}
-1 & 0 & 1\\
1 & -1 & 0\\
0 & 1 & -1
      \end{array}
\right)
\left(\begin{array}{ccc}
{l_1\over L_1} & 0 & 0\\
0 & {l_2\over L_2}  & 0\\
0 & 0 & {l_3\over L_3} 
      \end{array}
\right)
\left(
\begin{array}{ccc}
0 & \lambda_1 & \mu_1 \\ 
\mu_2 & 0 & \lambda_2 \\ 
\lambda_3 & \mu_3 & 0
\end{array} 
\right)
\left(\begin{array}{rrr}
-1 & 1 & 0 \\ 
0 & -1 & 1 \\ 
1 & 0 & -1
\end{array} \right)
\left(\begin{array}{c}
f_1 \\ 
f_2 \\ 
f_3
\end{array} \right)
   =
\left(\begin{array}{c}
A_1q_1 \\ 
A_2q_2 \\ 
A_3q_3
\end{array} \right).
\]
The diffusive term matrix 
\begin{eqnarray*}
&&
\left(\begin{array}{ccc}
-1 & 0 & 1\\
1 & -1 & 0\\
0 & 1 & -1
      \end{array}
\right)
\left(\begin{array}{ccc}
{l_1\over L_1} & 0 & 0\\
0 & {l_2\over L_2}  & 0\\
0 & 0 & {l_3\over L_3} 
      \end{array}
\right)
\left(
\begin{array}{ccc}
0 & \lambda_1 & \mu_1 \\ 
\mu_2 & 0 & \lambda_2 \\ 
\lambda_3 & \mu_3 & 0
\end{array} 
\right)
\left(\begin{array}{rrr}
-1 & 1 & 0 \\ 
0 & -1 & 1 \\ 
1 & 0 & -1
\end{array} \right)\\
&=&
\left(\begin{array}{ccc}
-{\lambda_3l_3\over L_3}-{\mu_1l_1\over L_1} & {\lambda_1l_1\over L_1}+{(\lambda_3-\mu_3)l_3\over L_3} 
          & -{(\lambda_1-\mu_1)l1\over L_1}+{\mu_3l_3\over L3} \\ 
{\mu_1l_1\over L_1}-{(\lambda_2-\mu_2)l_2\over L_2} & -{\lambda_1l1\over L_1} -{\mu_2l_2\over L_2}
          &  {(\lambda_1-\mu_1)l_1\over L_1}+{\lambda_2l_2\over L_2}\\ 
{(\lambda_2-\mu_2)l_2\over L_2} + {\lambda_3l_3\over L_3} & {\mu_2l_2\over L_2}-{(\lambda_3-\mu_3)l_3\over L_3} 
          & -{\lambda_2l_2\over L_2}-{\mu_3l_3\over L_3}
\end{array} 
\right)
\end{eqnarray*}
is actually symmetric (see Subsection \ref{subsubsec: diffusive term DEC FEML}).

\subsection{Local FEML-Discretized 2D anisotropic Poisson equation}

The diffusive elemental matrix in FEM (frequently called “stiffness matrix”) on an element $e$ is given by
\[K_e=\int B^tDBdA,\]
where $D$ is the matrix representing the anisotropic diffusion tensor $K$ in this paper, 
and the matrix $B$ is given explicitly by
\[
B=\left(
\begin{array}{ccc}
{\partial N_1\over\partial x} & {\partial N_2\over\partial x} & {\partial N_3\over\partial x} \\ 
{\partial N_1\over\partial y} & {\partial N_2\over\partial y} & {\partial N_3\over\partial y}
\end{array} 
\right)
=
{1\over 2A}
\left(
\begin{array}{ccc}
y_2-y_3 & y_3-y_1 & y_1-y_2 \\ 
x_3-x_2 & x_1-x_3 & x_2-x_1
\end{array} 
\right) . 
\]
Since the matrix $B$ is constant on an element of the mesh, the integral is easy to compute. Thus, the difussive matrix $K_e$ for a linear triangular element (FEML) is given by
\begin{eqnarray*}
K_e
   &=&
\int B^T DB dA \\
   &=&
  B^T D B A_e \\
   &=&
  {1\over 4A_e}
\left(
\begin{array}{ccc}
y_2-y_3 & x_3-x_2  \\ 
y_3-y_1 & x_1-x_3  \\
y_1-y_2 & x_2-x_1
\end{array} 
\right)
  \left(\begin{array}{cc}
  k_{11} & k_{12} \\ 
  k_{12} & k_{22}
  \end{array} \right)
\left(
\begin{array}{ccc}
y_2-y_3 & y_3-y_1 & y_1-y_2 \\ 
x_3-x_2 & x_1-x_3 & x_2-x_1
\end{array} 
\right)
\end{eqnarray*}
Now, let us consider the first diagonal entry of the local FEML anisotropic Poisson diffusive matrix $K_e$, 
\begin{eqnarray*}
  (K_e)_{11} 
    &=&
{1\over 4A}(k_{11}(y_2-y_3)^2 +(k_{12}+k_{12})(y_2-y_3)(x_3-x_2) + k_{22}(x_3-x_2)^2)\\
   &=& 
{1\over 4A}  (\begin{array}{cc}
  -(y_3-y_2), & x_3-x_2
  \end{array} )
  \left(\begin{array}{cc}
  k_{11} & k_{12} \\ 
  k_{12} & k_{22}
  \end{array} \right)
  \left(\begin{array}{c}
  -(y_3-y_2) \\ 
  x_3-x_2
  \end{array} \right)\\
   &=& 
{1\over 4A}  (\begin{array}{cc}
  x_3-x_2, & y_3-y_2
  \end{array} )
  \left(\begin{array}{cc}
  0 & 1 \\ 
  -1 & 0
  \end{array} \right)
  \left(\begin{array}{cc}
  k_{11} & k_{12} \\ 
  k_{12} & k_{22}
  \end{array} \right)
  \left(\begin{array}{cc}
  0 & -1 \\ 
  1 & 0
  \end{array} \right)
  \left(\begin{array}{c}
  x_3-x_2 \\ 
  y_3-y_2
  \end{array} \right)\\
   &=&
{1\over 4A}  (J(v_3-v_2))^T K J(v_3-v_2)\\
   &=&
{1\over 4A}  J(v_3-v_2)\cdot K (J(v_3-v_2)),
\end{eqnarray*}
where 
\begin{eqnarray*}
J
   &=&
  \left(\begin{array}{cc}
  0 & -1 \\ 
  1 & 0
  \end{array} \right)
\end{eqnarray*}
is the $90^\circ$ counter-clockwise rotation.
In this notation, the diffusive term in local FEML is given as follows
{\small
\[
{1\over 4A}
\left(\begin{array}{ccc}
J(v_3-v_2)\cdot K(J(v_3-v_2)) & J(v_3-v_2)\cdot K(J(v_1-v_3)) & J(v_3-v_2)\cdot K(J(v_2-v_1)) \\ 
J(v_1-v_3)\cdot K(J(v_3-v_2)) & J(v_1-v_3)\cdot K(J(v_1-v_3)) & J(v_1-v_3)\cdot K(J(v_2-v_1)) \\ 
J(v_2-v_1)\cdot K(J(v_3-v_2)) & J(v_2-v_1)\cdot K(J(v_1-v_3)) & J(v_2-v_1)\cdot K(J(v_2-v_1))
\end{array}\right) . \]
}

\subsection{Comparison between local DEC and FEML discretizations}
For the sake of brevity, 
we are only going to compare the entries of the first row and first column of each formulation.
Consider the various lengths, areas and angles given in the triangle of Figure \ref{fig: circumscribed triangle 01}.
    \begin{center}
    \includegraphics[width=0.33\textwidth]{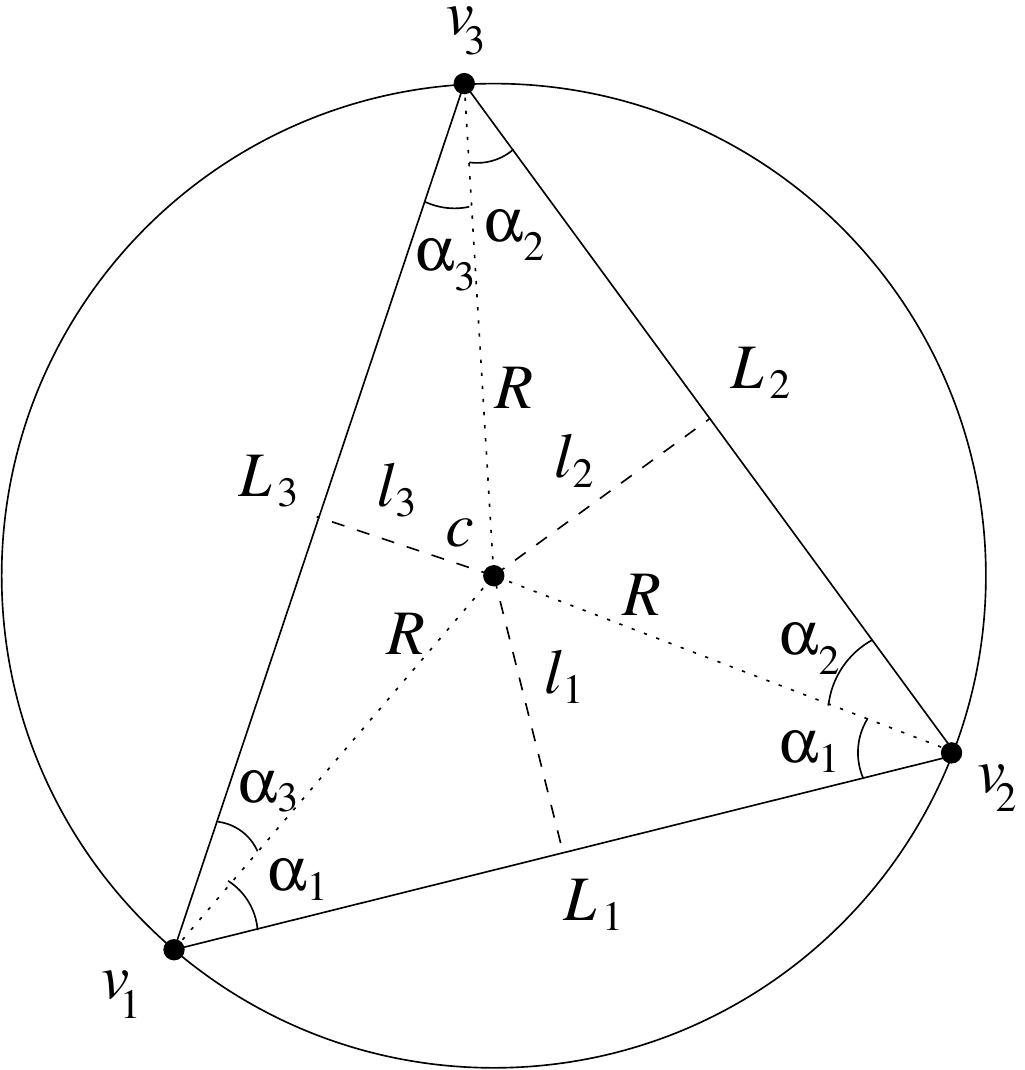}%
    \hspace{12mm}
    \includegraphics[width=0.33\textwidth]{WellCenteredTriangle02-eps-converted-to.pdf}%
    \captionof{figure}{Circumscribed triangle.}\label{fig: circumscribed triangle 01}
    \end{center}
We have the following:
\begin{eqnarray*}
 \pi &=& 2(\alpha_1+\alpha_2+\alpha_3),\\
 {2l_i\over L_i} &=& \tan(\alpha_i),\\
 {l_i\over R} &=& \sin(\alpha_i),\\
 {L_i\over 2R} &=& \cos(\alpha_i),\\
 A_1&=& {L_1l_1\over 4} + {L_3l_3\over 4},\\
 A_2&=& {L_1l_1\over 4} + {L_2l_2\over 4},\\
 A_3&=& {L_2l_2\over 4} + {L_3l_3\over 4}.
\end{eqnarray*}

\subsubsection{The diffusive term}\label{subsubsec: diffusive term DEC FEML}
We claim that
\[{J(v_3-v_2)\cdot K(J(v_3-v_2))\over 4A}
= 
-{\lambda_3l_3\over L_3}-{\mu_1l_1\over L_1}.
\]
Indeed,
\begin{eqnarray*}
\lambda_3
   &=&
  -{J(v_3-v_2)\cdot K( v_1-v_3)\over 2A},\\
\mu_1
   &=&
  {J(v_3-v_2)\cdot K( v_2-v_1)\over 2A}.
\end{eqnarray*}
Thus
\begin{eqnarray*}
-{\lambda_3l_3\over L_3}-{\mu_1l_1\over L_1}
   &=& 
   {J(v_3-v_2)\cdot K(v_1-v_3)\over 2A}{\tan(\alpha_3)\over 2}
   -{J(v_3-v_2)\cdot K( v_2-v_1)\over 2A}{\tan(\alpha_1)\over 2}\\
   &=& 
   {1\over 4A}J(v_3-v_2)\cdot K((v_1-v_3)\tan(\alpha_3)- (v_2-v_1)\tan(\alpha_1)).
\end{eqnarray*}
All we have to do is show that
\[(v_1-v_3)\tan(\alpha_3)- (v_2-v_1)\tan(\alpha_1) = J(v_3-v_2).\]
Note that, since $J(v_3-v_2)$ is orthogonal to $v_3-v_2$, $J(v_3-v_2)$ must be parallel to $c -  {v_2+v_3\over 2}$. Thus, 
\begin{equation}
J(v_3-v_2) = {L_2\over l_2}\left(c -  {v_2+v_3\over 2}\right). \label{eq: J(w2)}
\end{equation}
Now we are going to express $c$ in terms of $v_1,v_2,v_3$.
Let us consider 
\[c-v_1 = a (v_2-v_1) + b(v_3-v_1)\]
where $a,b$ are coefficients to be determined.
Taking inner products with $(v_2-v_1)$ and $(v_3-v_1)$ we get the two equations
\begin{eqnarray*}
R\cos(\alpha_1)
   &=&
  aL_1+bL_3\cos(\alpha_1+\alpha_3),\\ 
R\cos(\alpha_3)
   &=&
  aL_1\cos(\alpha_1+\alpha_3)+bL_3. 
\end{eqnarray*}
Solving for $a$ and $b$
\begin{eqnarray*}
a
   &=&
  {\sin(\alpha_3)\over 2\cos(\alpha_1)\sin(\alpha_1+\alpha_3)},\\
b
   &=&
  {\sin(\alpha_1)\over 2\cos(\alpha_3)\sin(\alpha_1+\alpha_3)}.
\end{eqnarray*}
Substituting all the relevant quantities in (\ref{eq: J(w2)}) we have, for instance, that the coefficient of $(v_2-v_1)$ is
\begin{eqnarray*}
2{\cos(\alpha_2)\over \sin(\alpha_2)}\left({\sin(\alpha_3)\over 2\cos(\alpha_1)\sin(\alpha_1+\alpha_3)}-{1\over 2}\right)
   &=&
{\cos(\alpha_2)\over \sin(\alpha_2)}\left({\sin(\alpha_3)-\cos(\alpha_1)\sin(\alpha_1+\alpha_3)\over \cos(\alpha_1)\sin(\alpha_1+\alpha_3)}\right)\\  
   &=&
{\cos(\alpha_2)\over \sin(\alpha_2)}\left({\sin(\alpha_3)-\cos(\alpha_1)(\sin(\alpha_1)\cos(\alpha_3)+\sin(\alpha_3)\cos(\alpha_1))\over \cos(\alpha_1)\sin(\pi/2-\alpha_2)}\right)\\  
   &=&
{\sin(\alpha_1)\over \sin(\alpha_2)}\left({\sin(\alpha_3)\sin(\alpha_1)-\cos(\alpha_1)\cos(\alpha_3)\over \cos(\alpha_1)}\right)\\  
   &=&
  \tan(\alpha_1){-\cos(\alpha_1+\alpha_3)\over\sin(\alpha_2)}\\ 
   &=&
  \tan(\alpha_1){-\cos(\pi/2-\alpha_2)\over\sin(\alpha_2)}\\ 
   &=&
  \tan(\alpha_1){-\sin(\alpha_2)\over\sin(\alpha_2)}\\ 
   &=&
  -\tan(\alpha_1), 
\end{eqnarray*}
and similarly for the coefficient of $(v_1-v_3)$.
The calculations for the remaining entries are similar.

Thus, the local DEC and FEML diffusive terms of the 2D anisotropic Poisson equation coincide.

\subsubsection{The source term}
As already observed in \cite{Esqueda1}, the right hand sides of the local DEC and FEML systems are different
\[\left(\begin{array}{c}
A_1q_1\\
A_2q_2\\
A_3q_3
      \end{array}
\right)
\not =
{A\over 3}\left(\begin{array}{l}
q_1\\
q_2\\
q_3
        \end{array}
\right).\]
While FEML uses a barycentric subdivision to calculate the areas associated to each node/vertex,
DEC uses a circumcentric subdivision. Eventually, this leads the DEC discretization to a better approximation of the solution
(on coarse meshes).

\section{Some remarks about DEC quantities}\label{sec: extension}

\subsection{The discrete Hodge star quantities revisited}   \label{subsec: determinants} 
The numbers appearing in the local DEC matrices can be expressed both in terms of determinants and in terms of trigonometric functions.
More precisely, 
\begin{eqnarray*}
 A_1
&=& {1\over 4}\left[\det\left(\begin{array}{ccc}
x_1 & y_1 & 1\\
x_c & y_c & 1\\
x_2 & y_2 & 1
                            \end{array}
\right)
+\det\left(\begin{array}{ccc}
x_3 & y_3 & 1\\
x_c & y_c & 1\\
x_1 & y_1 & 1
           \end{array}
\right)
\right]
 \quad =\quad {R^2\over 4}(\sin(2\alpha_1)+\sin(2\alpha_3)),\\
 A_2
 &=& {1\over 4}\left[\det\left(\begin{array}{ccc}
x_1 & y_1 & 1\\
x_c & y_c & 1\\
x_2 & y_2 & 1
                            \end{array}
\right)
+\det\left(\begin{array}{ccc}
x_2 & y_2 & 1\\
x_c & y_c & 1\\
x_3 & y_3 & 1
           \end{array}
\right)
\right]
\quad = \quad {R^2\over 4}(\sin(2\alpha_1)+\sin(2\alpha_2)),\\
 A_3
 &=& {1\over 4}\left[\det\left(\begin{array}{ccc}
x_2 & y_2 & 1\\
x_c & y_c & 1\\
x_3 & y_3 & 1
                            \end{array}
\right)
+\det\left(\begin{array}{ccc}
x_3 & y_3 & 1\\
x_c & y_c & 1\\
x_1 & y_1 & 1
           \end{array}
\right)
\right]
\quad = \quad {R^2\over 4}(\sin(2\alpha_2)+\sin(2\alpha_3)),\\
{l_1\over L_1}
 &=&
{1\over L_1^2}
\det\left(\begin{array}{ccc}
x_1 & y_1 & 1\\
x_c & y_c & 1\\
x_2 & y_2 & 1
                            \end{array}
\right)
\quad = \quad   {\tan(\alpha_1)\over 2},\\
{l_2\over L_2}
 &=&
{1\over L_2^2}
\det\left(\begin{array}{ccc}
x_2 & y_2 & 1\\
x_c & y_c & 1\\
x_3 & y_3 & 1
                            \end{array}
\right)
 \quad = \quad  {\tan(\alpha_2)\over 2},\\
{l_3\over L_3}
 &=&
{1\over L_3^2}
\det\left(\begin{array}{ccc}
x_3 & y_3 & 1\\
x_c & y_c & 1\\
x_1 & y_1 & 1
                            \end{array}
\right)
  \quad = \quad {\tan(\alpha_3)\over 2}.
\end{eqnarray*}
These expressions are valid regardless of the location of the circumcenter and can, indeed, take negative values.
he angles that are measured in the scheme can be negative as in the obtuse triangle of Figure \ref{fig: negative angles}
    \begin{center}
    \includegraphics[width=0.3\textwidth]{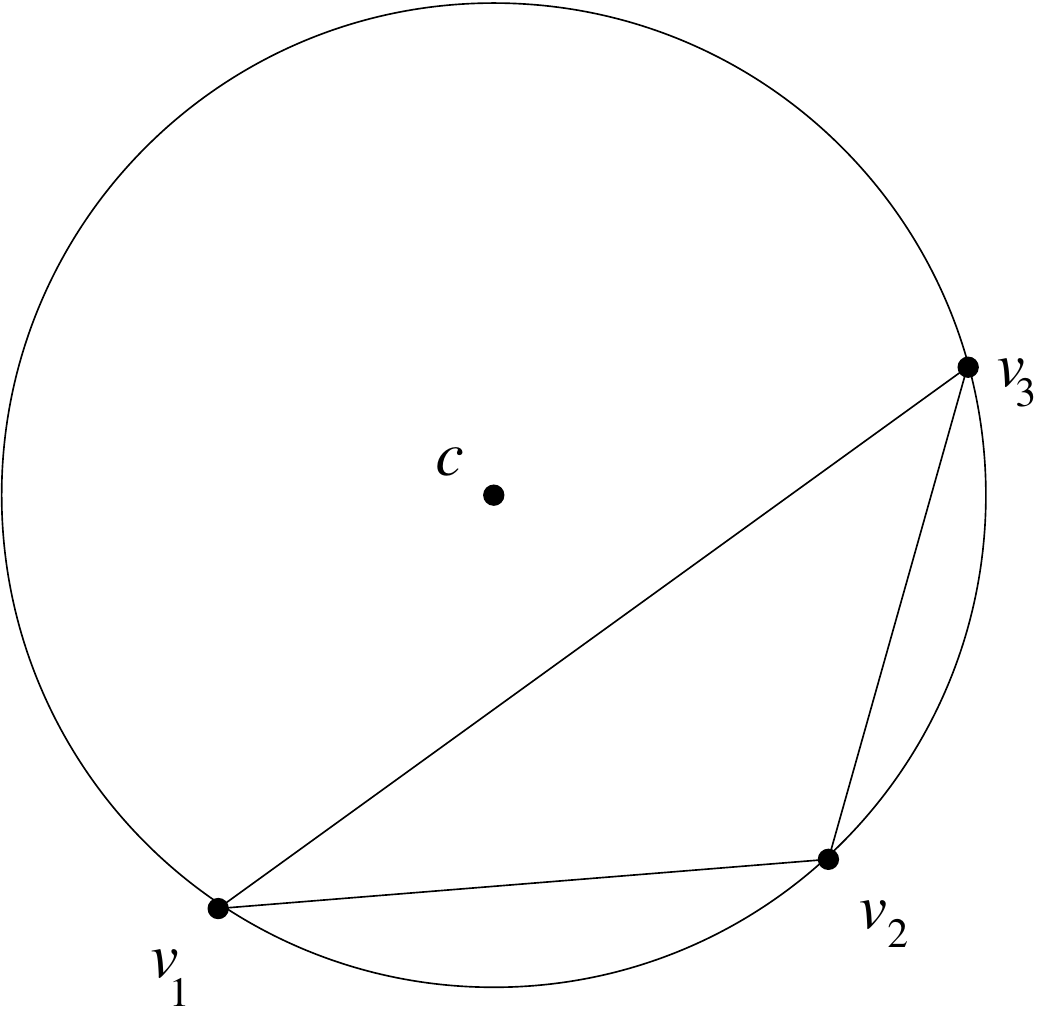}%
   \hspace{8mm}
    \includegraphics[width=0.3\textwidth]{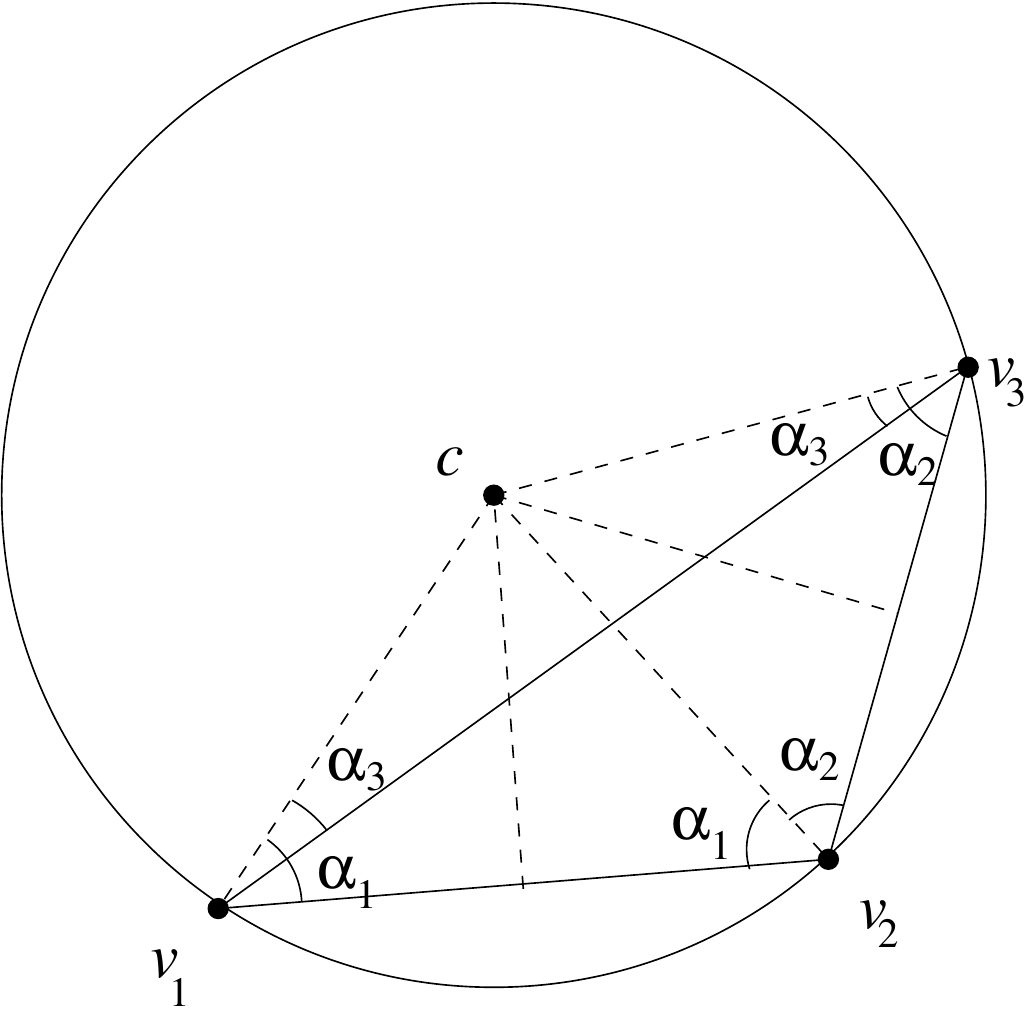}%
    \captionof{figure}{Negative (exterior) angles measured in an obtuse triangle.}\label{fig: negative angles}
    \end{center}
and some quantities can be zero or negative.
For instance, if 
\[\alpha_2={\pi\over 2}-2\alpha_1,\]
hen
\[A_1=0.\]

\subsection{Area weights assigned to vertices}
In order to understand how local DEC assigns area weights to vertices differently from FEML, let us
consider the obtuse triangle shown in Figure \ref{fig: negative angles}
Let $p_1,p_2,p_3$ be the middle points of the segments
$[v_2,v_3],[v_2,v_3],[v_3,v_1]$ respectively.    
As shown in Figure \ref{fig: obtuse triangle 01}, the triangle $[v_1,p_3,c]$ lies completely outside of the triangle $[v_1,v_2,v_3]$. Geometrically, 
this implies that its area must be assigned a negative sign, which is confirmed 
by the determinant formulas of Subsection \ref{subsec: determinants}. On the other hand, the triangle $[v_1,p_1,c]$ will have positive area. 
Thus, their sum gives us the area $A_1$ in Figure \ref{fig: obtuse triangle 01}.
    \begin{center}
    \includegraphics[width=0.28\textwidth]{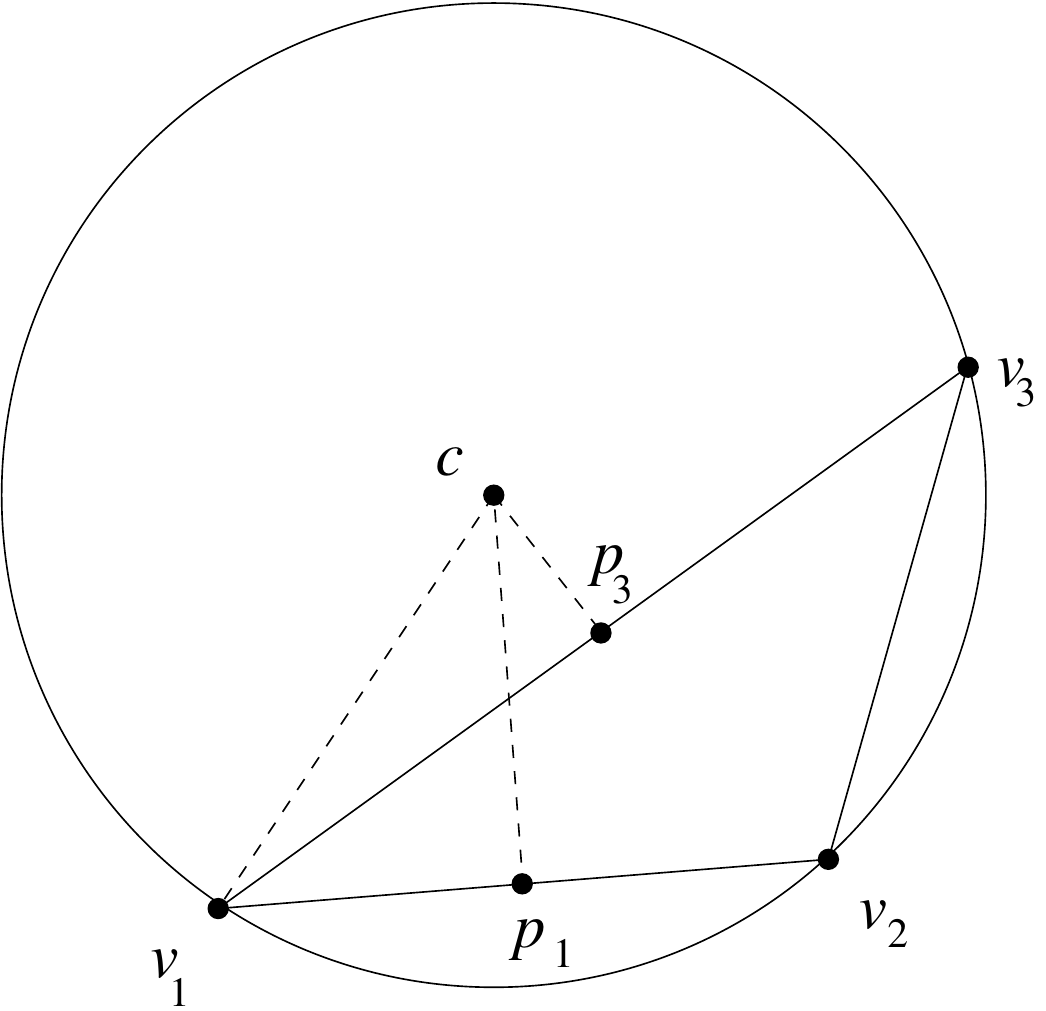}%
    \hspace{3mm}
    \includegraphics[width=0.28\textwidth]{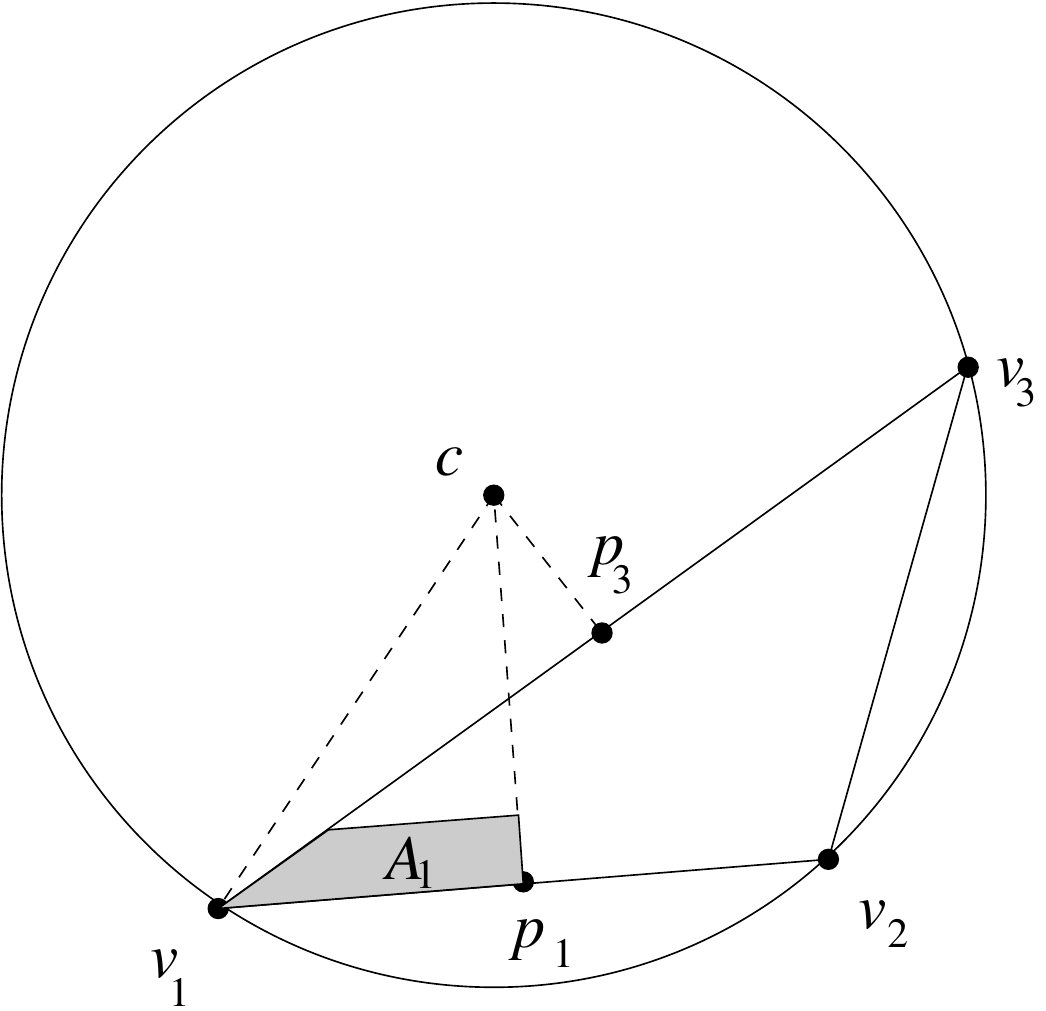}%
    \hspace{3mm}
    \includegraphics[width=0.28\textwidth]{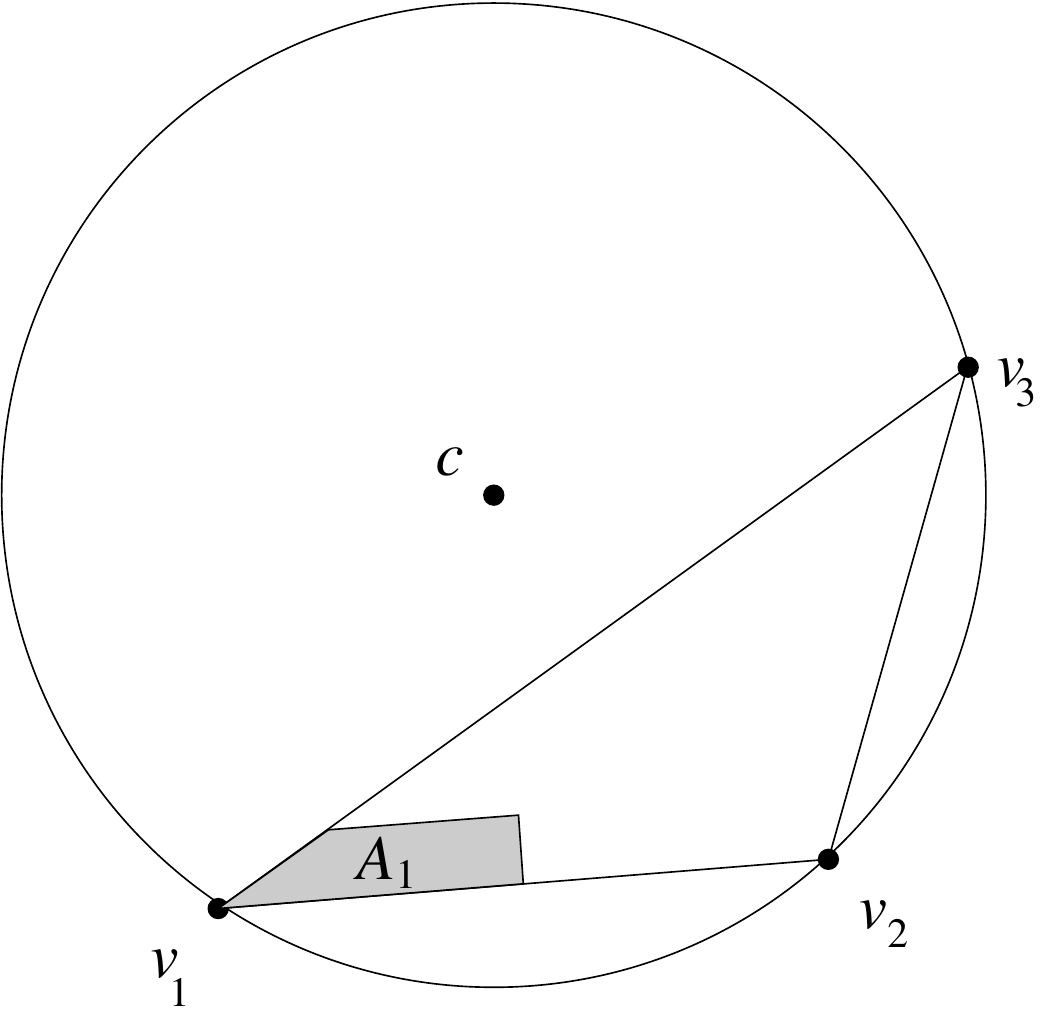}%
    \captionof{figure}{Area weight assigned to $v_1$.}\label{fig: obtuse triangle 01}
    \end{center}
The area $A_3$ is computed similarly, where the triangle $[p_3,v_3,c]$ is assigned negative area (see Figure \ref{fig: obtuse triangle 02}).
    \begin{center}    
    \includegraphics[width=0.28\textwidth]{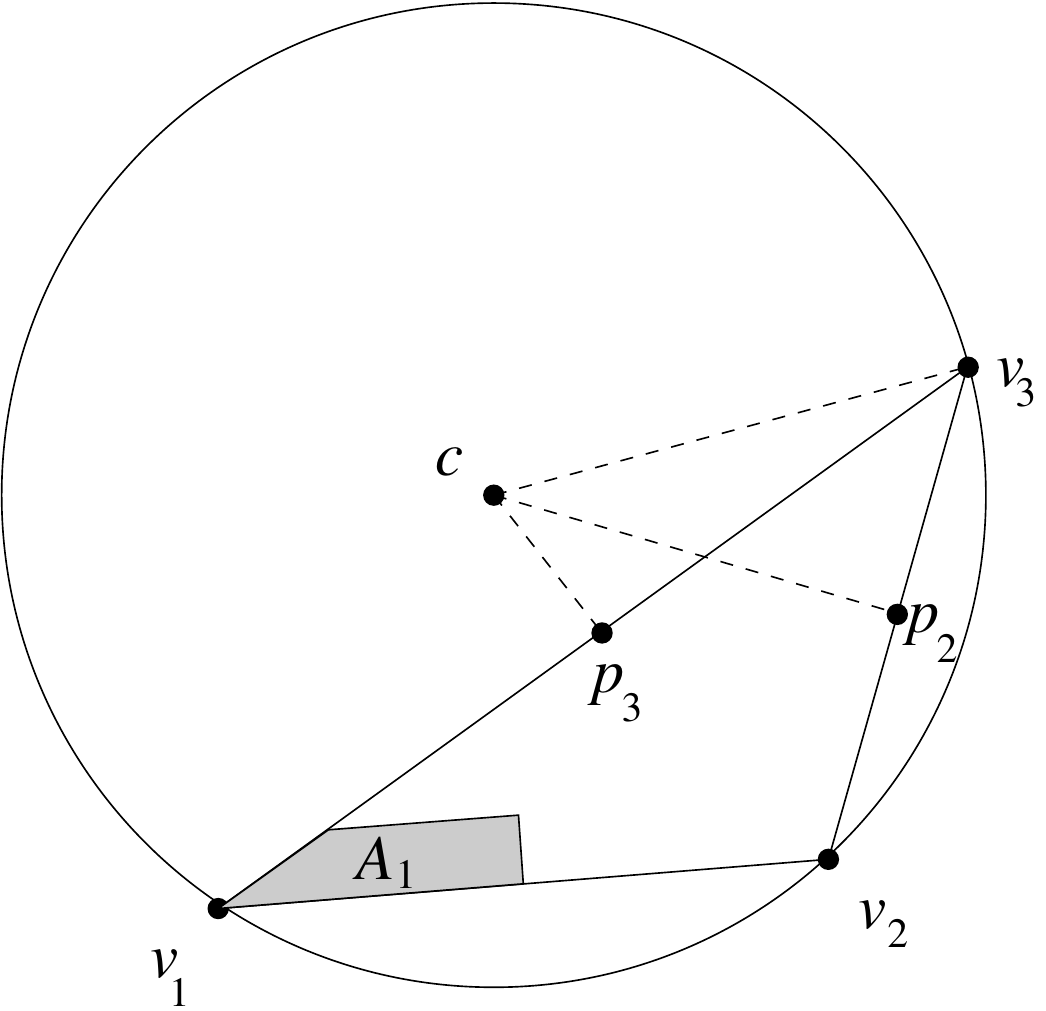}%
    \hspace{3mm}
    \includegraphics[width=0.28\textwidth]{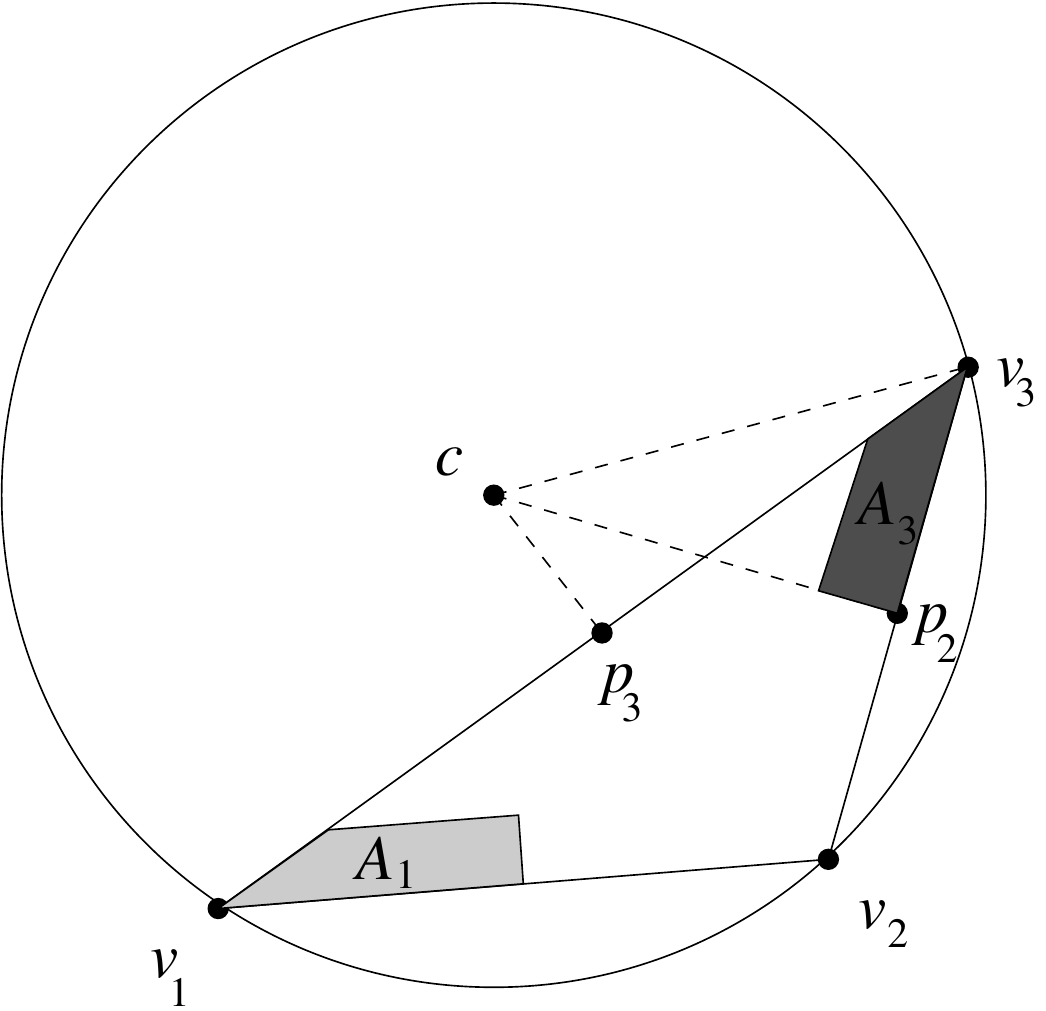}%
    \hspace{3mm}
    \includegraphics[width=0.28\textwidth]{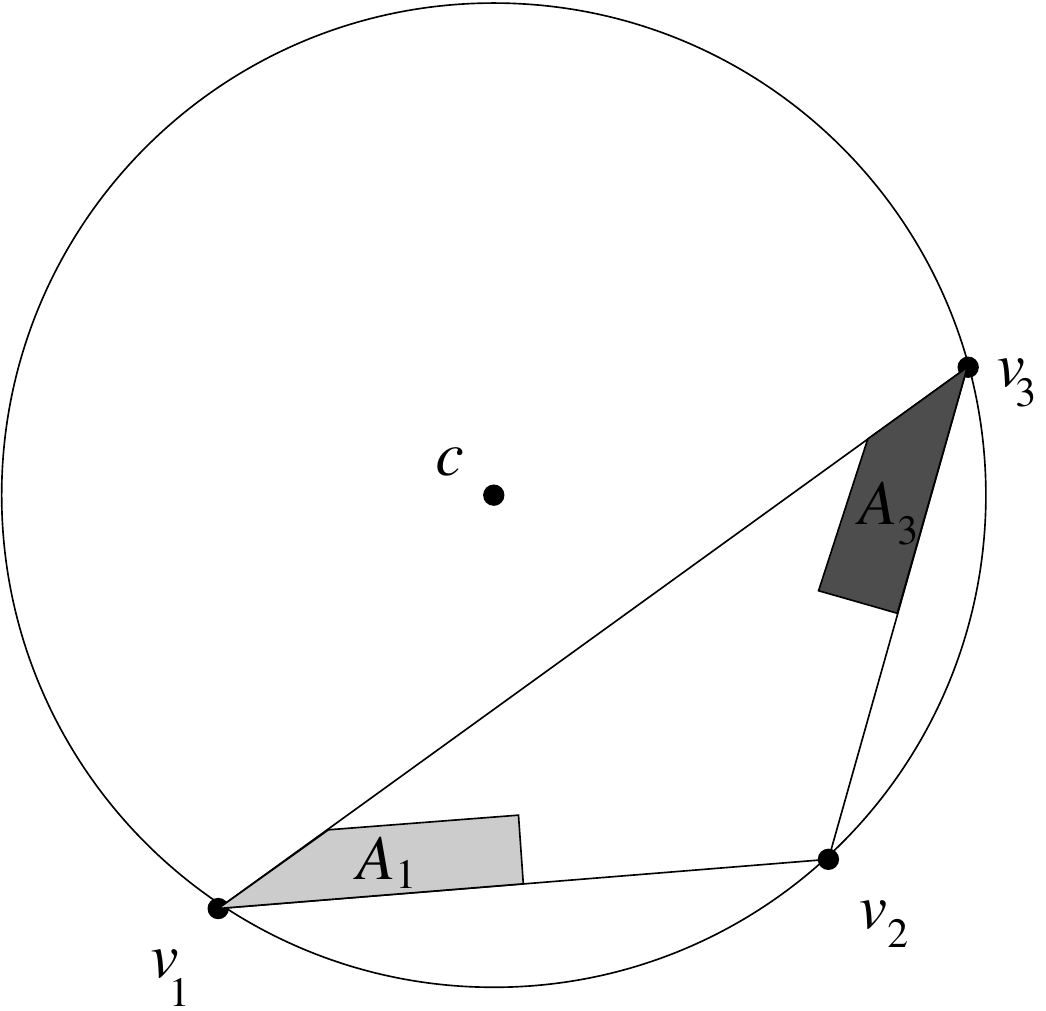}%
    \captionof{figure}{Area weight assigned to $v_3$.}\label{fig: obtuse triangle 02}
    \end{center}
Note that for $A_2$, the two triangles $[p_1,v_3,c]$ and $[v_2,p_2,c]$ both have positive areas (see Figure \ref{fig: obtuse triangle 03}). 
%The phenomenon of one node taking the lion's share of the area is typical in obtuse triangles.
    \begin{center}
    \includegraphics[width=0.28\textwidth]{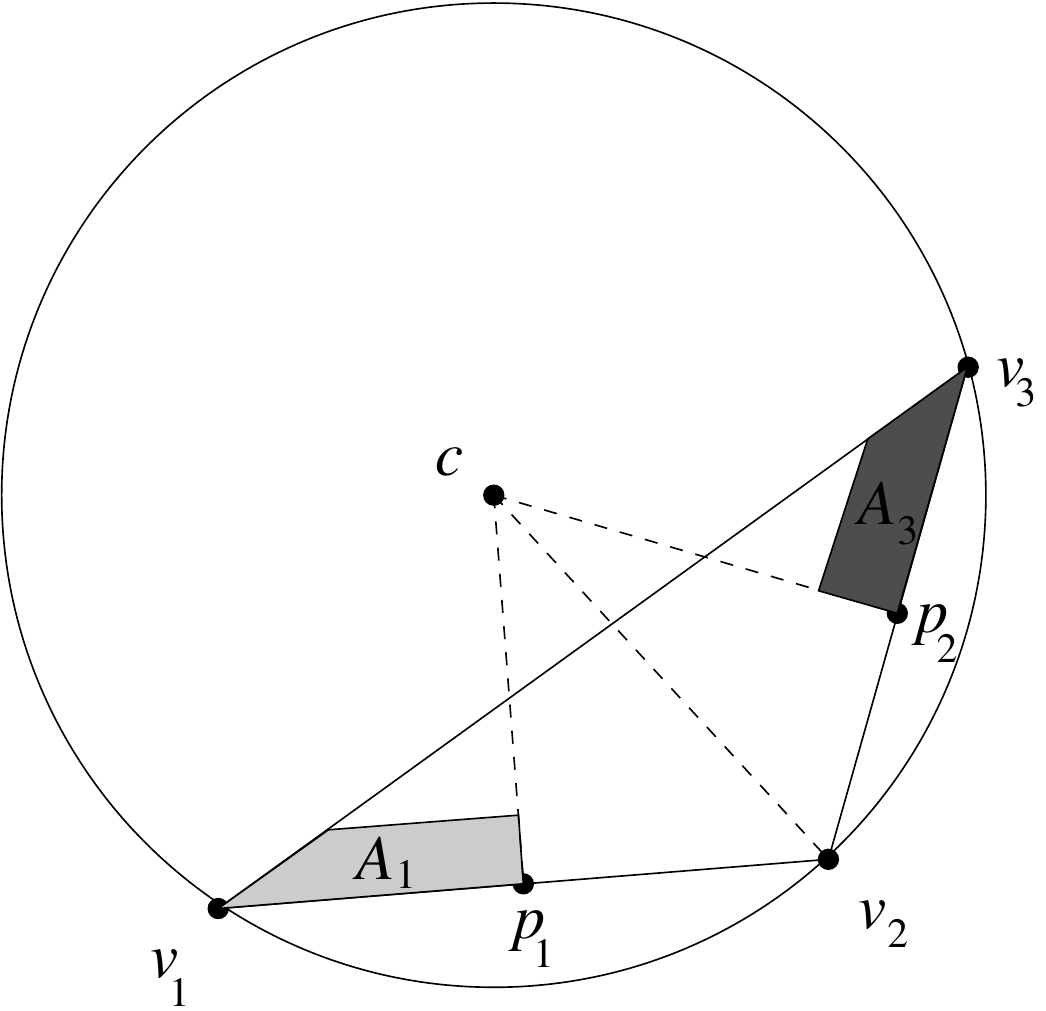}%
    \hspace{2mm}
    \includegraphics[width=0.28\textwidth]{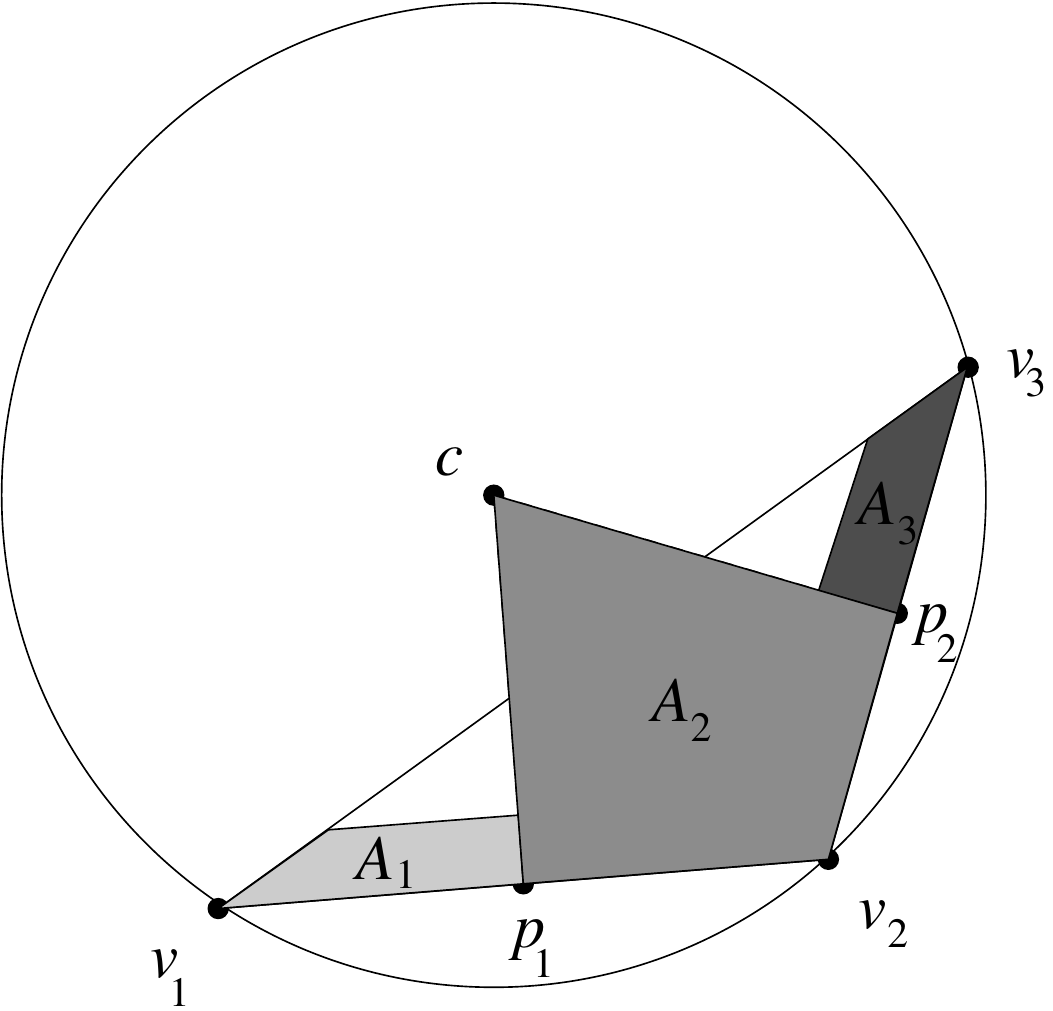}%
    \hspace{2mm}
    \includegraphics[width=0.28\textwidth]{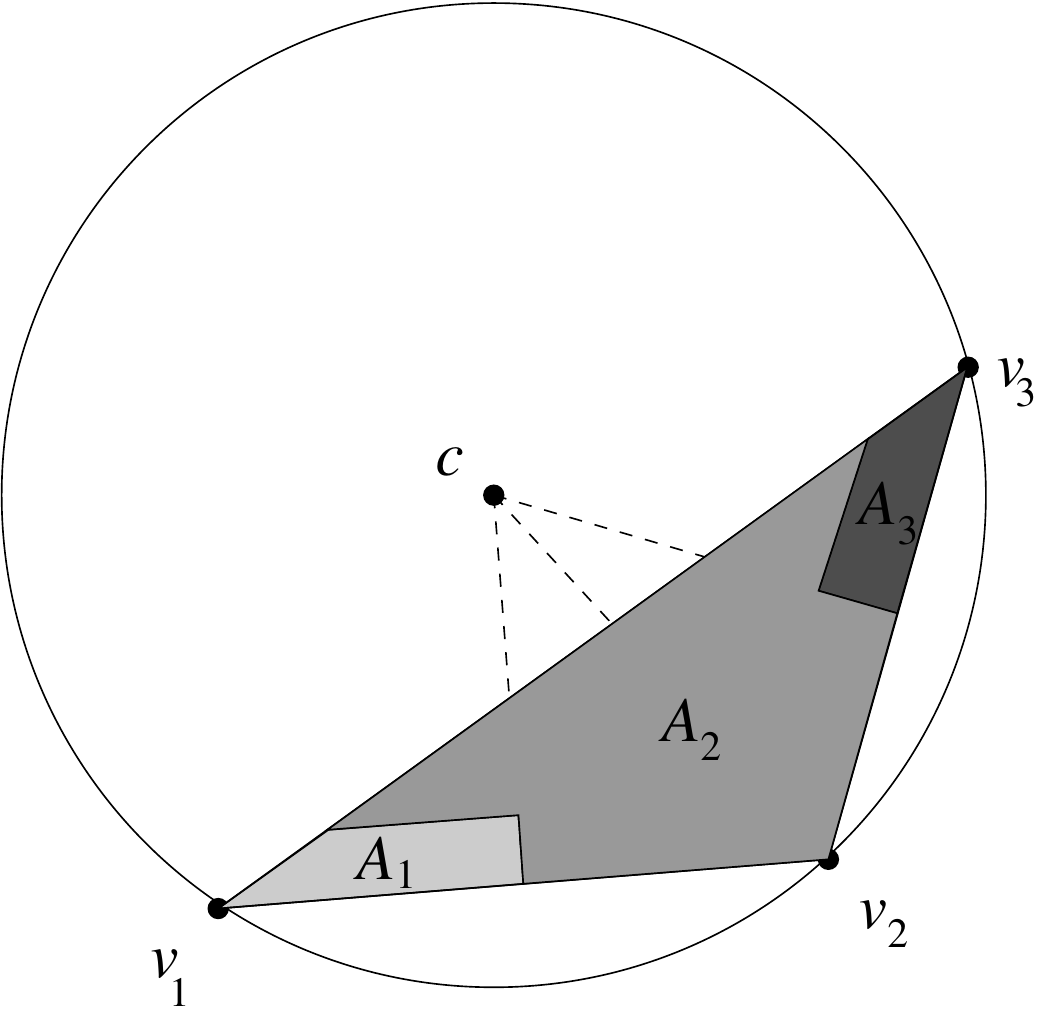}%
    \captionof{figure}{Area weight assigned to $v_2$.}\label{fig: obtuse triangle 03}
    \end{center}

\section{Numerical Examples}\label{sec: examples}

In this section, we present three examples in order to illustrate the performance of DEC resulting from the local formulation and its implementation. 
In all cases, we solve the anisotropic Poisson equation.
The FEML methodology that we have used in the comparison can be consulted \cite{Onate,Zienkiewicz1,Botello}.

\subsection{First example: Heterogeneity}
\begin{figure}[h]
	\centering
	\includegraphics[width=0.38\textwidth]{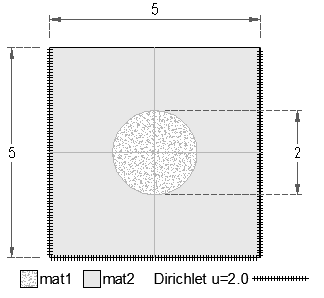}%
	\captionof{figure}{Square and inner circle with different conditions.}\label{fig: square example geometry}
\end{figure}

\begin{figure}[h]
	\centering
	\subfigure[]{\includegraphics[width=0.25\textwidth]{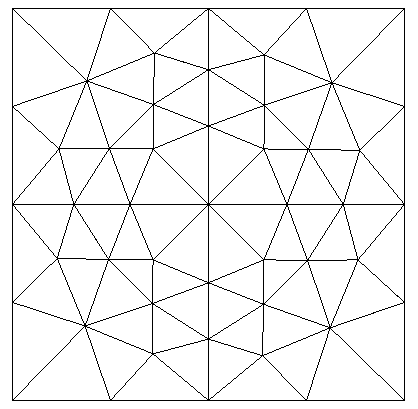}}\hspace{5mm}
	\subfigure[]{\includegraphics[width=0.25\textwidth]{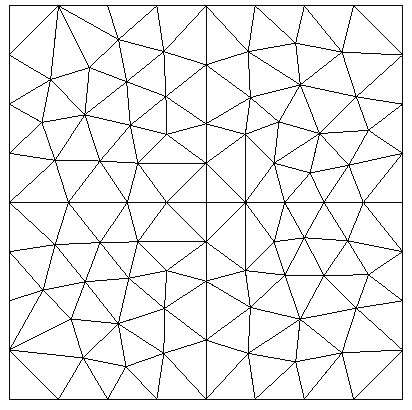}}\hspace{5mm}
	\subfigure[]{\includegraphics[width=0.25\textwidth]{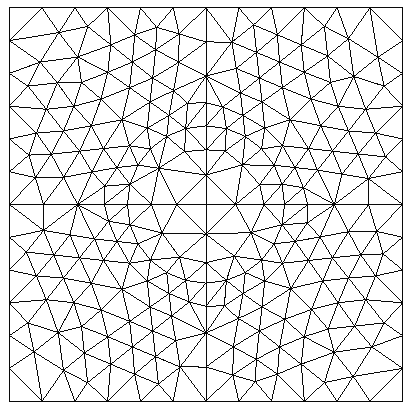}}
	%\vspace{2mm}
	\subfigure[]{\includegraphics[width=0.25\textwidth]{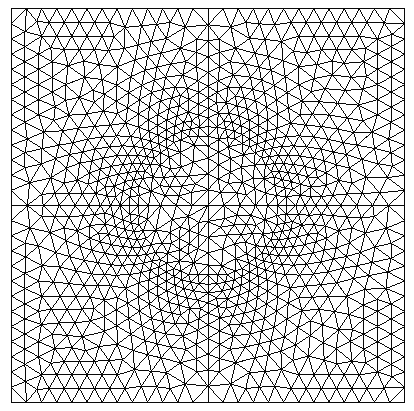}}\hspace{5mm}
	\subfigure[]{\includegraphics[width=0.25\textwidth]{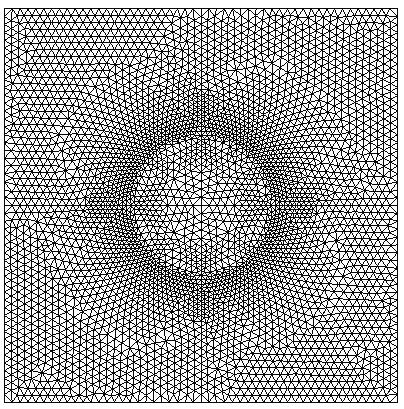}}\hspace{5mm}
	\subfigure[]{\includegraphics[width=0.25\textwidth]{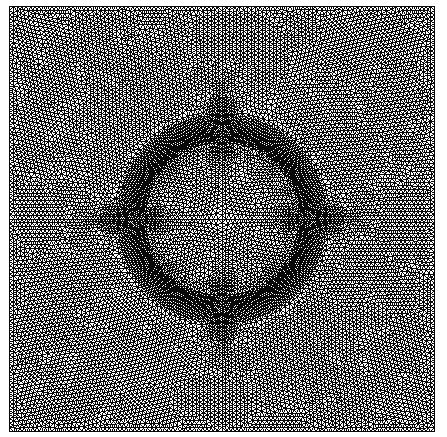}}
	%\vspace{2mm}
	\captionof{figure}{Six of the meshes used in the first example.}\label{fig: first example meshes}
\end{figure}

This example is intended to highlight how Local DEC deals effectively with heterogeneous materials. 
Consider the region in the plane given in Figure \ref{fig: square example geometry}.
\begin{itemize}
	\item The difussion constant for the region labelled mat1 is $k=12$ and its source term is $q=20$.
	\item The difussion constant for the region labelled mat2 is $k=6$ and its source term is $q=5$.
\end{itemize}
The meshes used in this example are shown in Figure \ref{fig: first example meshes} and vary from coarse to very fine.

\begin{figure}[h]
	\centering
	\subfigure[Contour Fill of Temperatures]{\includegraphics[width=0.4\textwidth]{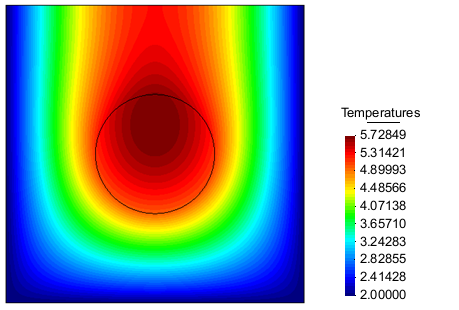}}\hspace{3mm}
	\subfigure[Contour Fill of Flux vectors on Elems]{\includegraphics[width=0.4\textwidth]{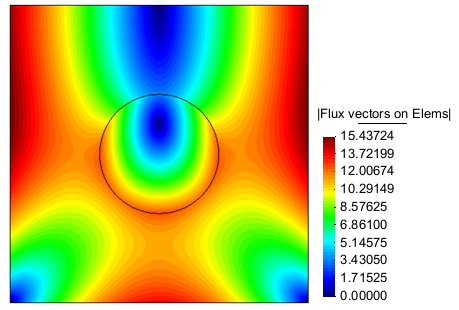}}
	\captionof{figure}{Temperature and flux-magnitude distribution fields of the first example.}\label{fig: square example temperature field}
\end{figure}

The numerical results for the maximum temperature value are exemplified in Table \ref{table: square example numerical results}.
\begin{center}
	\begin{tabular}{|c|c|c|c|c|c|c|}\hline
		\multirow{2}{*}{\rm Mesh}&\multirow{2}{*}{	\# {\rm nodes}}&\multirow{2}{*}{\# {\rm elements}}	&	\multicolumn{2}{c|}{\mbox{Max. Temp. Value}}&	\multicolumn{2}{c|}{\mbox{Max. Flux Magnitude}} \\
		&&&{\rm DEC} 	&{\rm FEML}&{\rm DEC} 	&{\rm FEML}	 \\\hline
		Figure \ref{fig: first example meshes}(a) &
		49 &
		80 &
		5.51836 &
		5.53345 &
		13.837 &
		13.453 \\
		Figure \ref{fig: first example meshes}(b) &
		98 &
		162 &
		5.65826 &
		5.66648 &
		14.137 &
		14.024 \\
		Figure \ref{fig: first example meshes}(c) &
		258 &
		466 &
		5.70585 &
		5.71709 &
		14.858 &
		14.770 \\
		Figure \ref{fig: first example meshes}(d) &
		1,010 &
		1,914 &
		5.72103 &
		5.72280 &
		15.008 &
		15.006 \\
		Figure \ref{fig: first example meshes}(e) &
		3,813 &
		7,424 &
		5.72725 &
		5.72725 &
		15.229 &
		15.228 \\
		Figure \ref{fig: first example meshes}(f) &
		13,911 &
		27,420 &
		5.72821 &
		5.72826 &
		15.342 &
		15.337 \\
		&
		50,950 &
		101,098 &
		5.72841 &
		5.72842 &
		15.395 &
		15.396 \\
		&
		135,519 &
		269,700 &
		5.72845 &
		5.72845 &
		15.420 &
		15.417 \\
		&
		298,299 &
		594,596 &
		5.72848 &
		5.72848 &
		15.430 &
		15.429 \\
		&
		600,594 &
		1,198,330 &
		5.72848 &
		5.72848 &
		15.433 &
		15.433 \\
		&
		1,175,238 &
		2,346,474 &
		5.72849 &
		5.72849 &
		15.43724 &
		15.43724 \\
		\hline
	\end{tabular}
	\captionof{table}{Numerical simulation results of the first example.}\label{table: square example numerical results}  
\end{center}
The temperature and flux-magnitude distribution fields are shown in Figure \ref{fig: square example temperature field}.    

\begin{figure}
	\centering
	\subfigure[]{\includegraphics[width=0.41\textwidth]{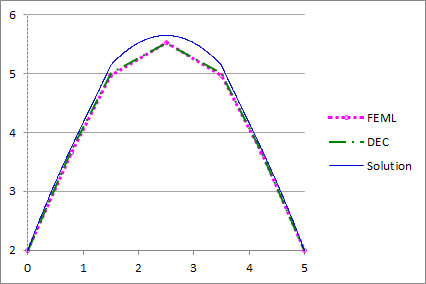}}\hspace{2mm}
	\subfigure[]{\includegraphics[width=0.41\textwidth]{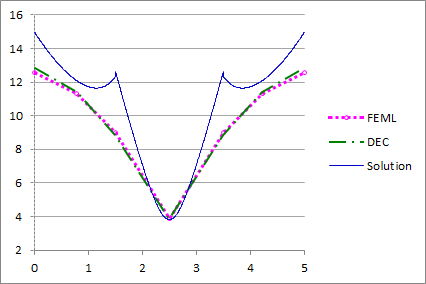}}
	\subfigure[]{\includegraphics[width=0.41\textwidth]{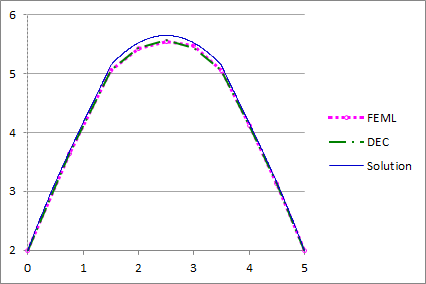}}\hspace{2mm}
	\subfigure[]{\includegraphics[width=0.41\textwidth]{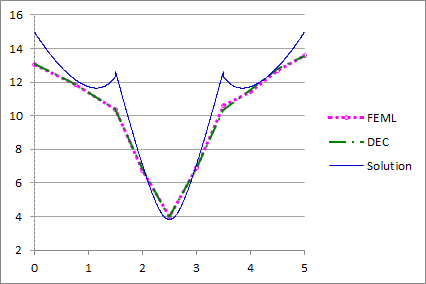}}
	\subfigure[]{\includegraphics[width=0.41\textwidth]{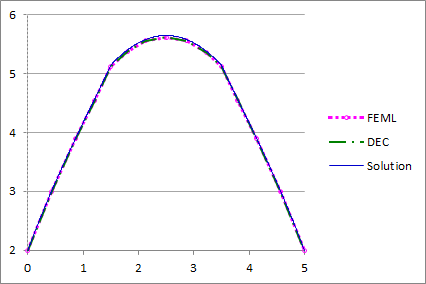}}\hspace{2mm}
	\subfigure[]{\includegraphics[width=0.41\textwidth]{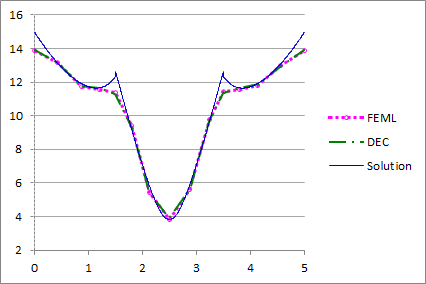}}
	\captionof{figure}{Temperature and Flux magnitude graphs of the first example along a cross-section of the domain for different meshes.}\label{fig: first example diametral graphs}
\end{figure}

Figure \ref{fig: first example diametral graphs} shows the graphs of the temperature and the flux-magnitude along a horizontal line crossing the inner circle for the first two meshes.

\newpage

\subsection{Second example: Anisotropy}
\begin{figure}[h]
	\centering
	\includegraphics[width=0.3\textwidth]{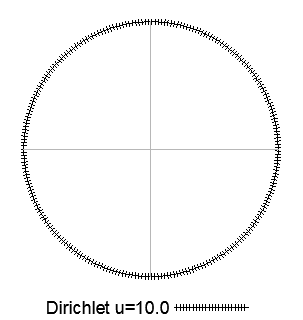}
	\captionof{figure}{Disk of radius one.}\label{fig: second example geometry}
\end{figure}

\begin{figure}[h]
	\centering
	\subfigure[]{\includegraphics[width=0.25\textwidth]{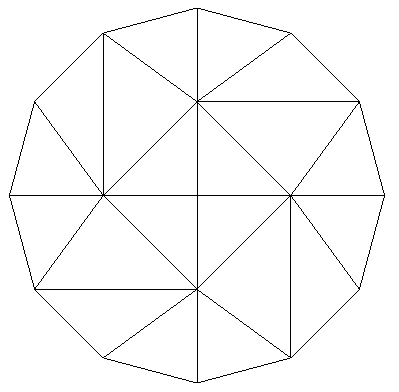}}\hspace{4mm}
	\subfigure[]{\includegraphics[width=0.25\textwidth]{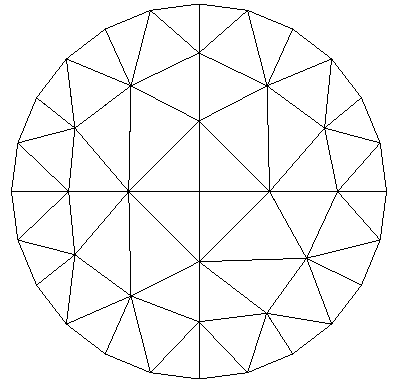}}\hspace{4mm}
	\subfigure[]{\includegraphics[width=0.25\textwidth]{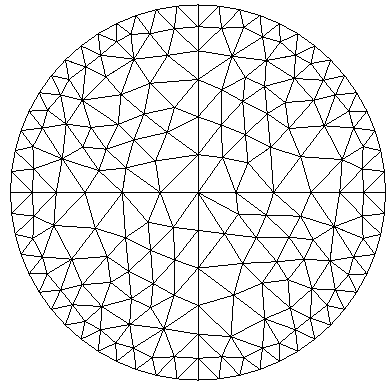}}
	\subfigure[]{\includegraphics[width=0.25\textwidth]{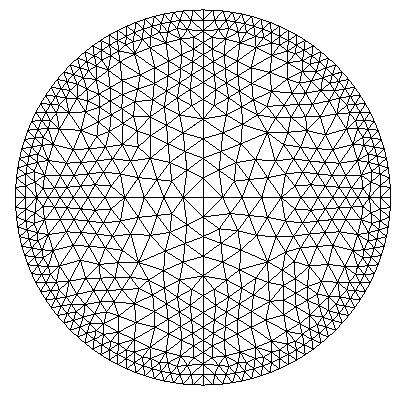}}\hspace{4mm}
	\subfigure[]{\includegraphics[width=0.25\textwidth]{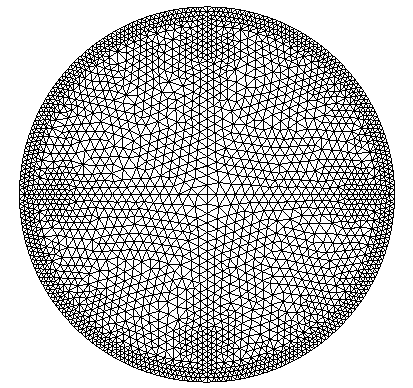}}\hspace{4mm}
	\subfigure[]{\includegraphics[width=0.25\textwidth]{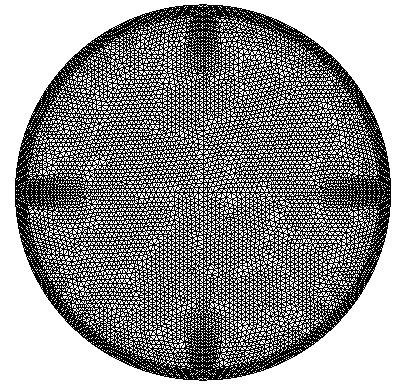}}
	%\vspace{2mm}
	\captionof{figure}{Six firsts meshes used for unit disk.}\label{fig: second example meshes}
\end{figure}

\begin{figure}[h]
	\centering
	\subfigure[Contour Fill of Temperatures]{\includegraphics[width=0.48\textwidth]{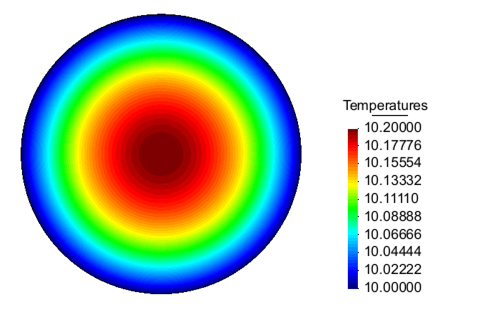}}\hspace{2mm}
	\subfigure[Contour Fill of Flux vectors on Elems]{\includegraphics[width=0.48\textwidth]{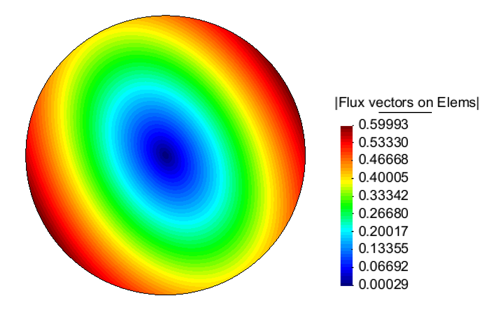}}
	\captionof{figure}{Temperature distribution and Flux magnitude fields for the finest mesh of the second example.}\label{fig: second example temperature field}
\end{figure}

Let us solve the Poisson equation in a circle of radius one centered at the origin $(0,0)$ under the following conditions (see Figure \ref{fig: second example geometry}): 
\begin{itemize}
	\item heat anisotropic diffusion constants $K_x = 1.5, K_y=1.0$;
	\item material angle $30^\circ$;
	\item source term $q= 1$;
	\item Dirichlet boundary condition $u=10$.
\end{itemize}
The meshes used in this example are shown in Figure \ref{fig: second example meshes} and vary from very coarse to very fine.
The numerical results for the maximum temperature value ($u(0,0)=10.2$) are exemplified in Table \ref{table: numerical results} where 
a comparison with the Finite Element Method with linear interpolation functions (FEML) is also shown. 
\begin{center}
	\begin{tabular}{|c|c|c|c|c|c|c|}\hline
		\multirow{2}{*}{\rm Mesh}&\multirow{2}{*}{	\# {\rm nodes}}&\multirow{2}{*}{\# {\rm elements}}	&	\multicolumn{2}{c|}{\mbox{Temp. Value at $(0,0)$}}&	\multicolumn{2}{c|}{\mbox{Flux Magnitude at $(-1,0)$}} \\
		&&&{\rm DEC} 	&{\rm FEML}&{\rm DEC} 	&{\rm FEML}	 \\\hline
		{\rm Figure \,\,\,\ref{fig: second example meshes}(a)}&	17	&	20		&	10.20014&			10.19002		& 0.42133 & 0.43865\\
		{\rm Figure \,\,\,\ref{fig: second example meshes}(b)}&	41	&	56		&	10.20007&			10.19678		& 0.48544 & 0.49387\\
		{\rm Figure \,\,\,\ref{fig: second example meshes}(c)}&	201	&	344		&	10.20012&			10.20158		& 0.52470 & 0.52428\\
		{\rm Figure \,\,\,\ref{fig: second example meshes}(d)}&	713	&	1304	&	10.20000&			10.19969		& 0.54143 & 0.54224\\
		{\rm Figure \,\,\,\ref{fig: second example meshes}(e)}&	2455&	4660	&	10.20000&			10.19990		& 0.54971 & 0.55138\\
		{\rm Figure \,\,\,\ref{fig: second example meshes}(f)}& 	8180&	15862	&	10.20000&			10.20002		& 0.55326 & 0.55409\\
		& 20016	&	39198	&	10.20000&			10.19999		& 0.55470 & 0.55520\\
		& 42306	&	83362	&	10.20000&			10.20000		& 0.55540 & 0.55572\\\hline
	\end{tabular}    
	\captionof{table}{Temperature value at the point $(0,0)$ and Flux magnitude value at the point $(-1,0)$ of the numerical simulations for the second example.}\label{table: numerical results} 
\end{center}
The temperature distribution and Flux magnitude fields for the finest mesh are shown in Figure \ref{fig: second example temperature field}.

\begin{figure}
	\centering
	\subfigure[]{\includegraphics[width=0.45\textwidth]{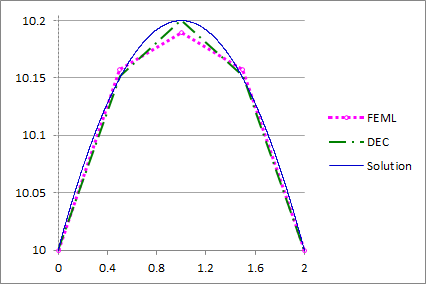}}\hspace{2mm}
	\subfigure[]{\includegraphics[width=0.45\textwidth]{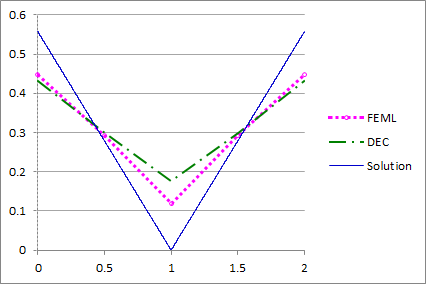}}\vspace{1mm}
	\subfigure[]{\includegraphics[width=0.45\textwidth]{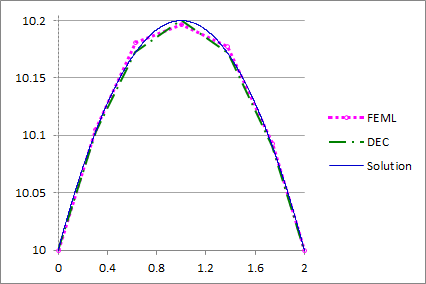}}\hspace{2mm}
	\subfigure[]{\includegraphics[width=0.45\textwidth]{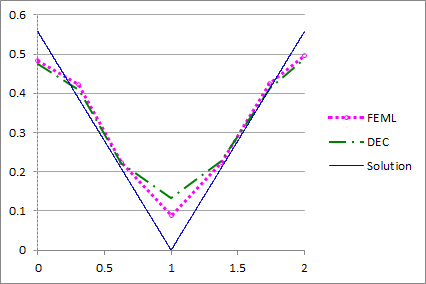}}\vspace{1mm}
	\subfigure[]{\includegraphics[width=0.45\textwidth]{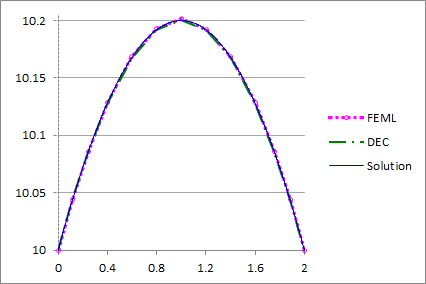}}\hspace{2mm}
	\subfigure[]{\includegraphics[width=0.45\textwidth]{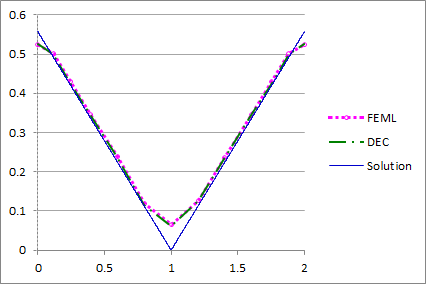}}\vspace{1mm}
	\captionof{figure}{Temperature and Flux magnitude graphs of the second example along a diameter of the circle for different meshes: 
		mesh in Figure \ref{fig: second example meshes}(a), a-Temperature, b-Flux; 
		mesh in Figure \ref{fig: second example meshes}(b), c-Temperature, d-Flux; 
		mesh in Figure \ref{fig: second example meshes}(c), e-Temperature, f-Flux; 
	}\label{fig: second example diametral graphs}
\end{figure}

Figures \ref{fig: second example diametral graphs}(a), \ref{fig: second example diametral graphs}(b) and \ref{fig: second example diametral graphs}(c) 
show the graphs of the temperature and flux magnitude values along a diameter of the circle for the different meshes of Figures \ref{fig: second example meshes}(a), \ref{fig: second example meshes}(b)
and \ref{fig: second example meshes}(c) respectively.

\newpage

\subsection{Third example: Heterogeneity and anisotropy}
\begin{figure}[h]
	\centering
	\includegraphics[width=0.4\textwidth]{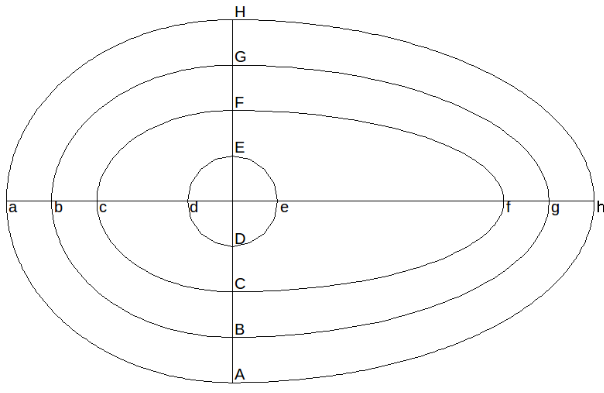}
	\captionof{figure}{Egg-like domain with different materials.}\label{fig: third example geometry}
\end{figure}

\begin{figure}[h]
	\centering
	\subfigure[]{\includegraphics[width=0.35\textwidth]{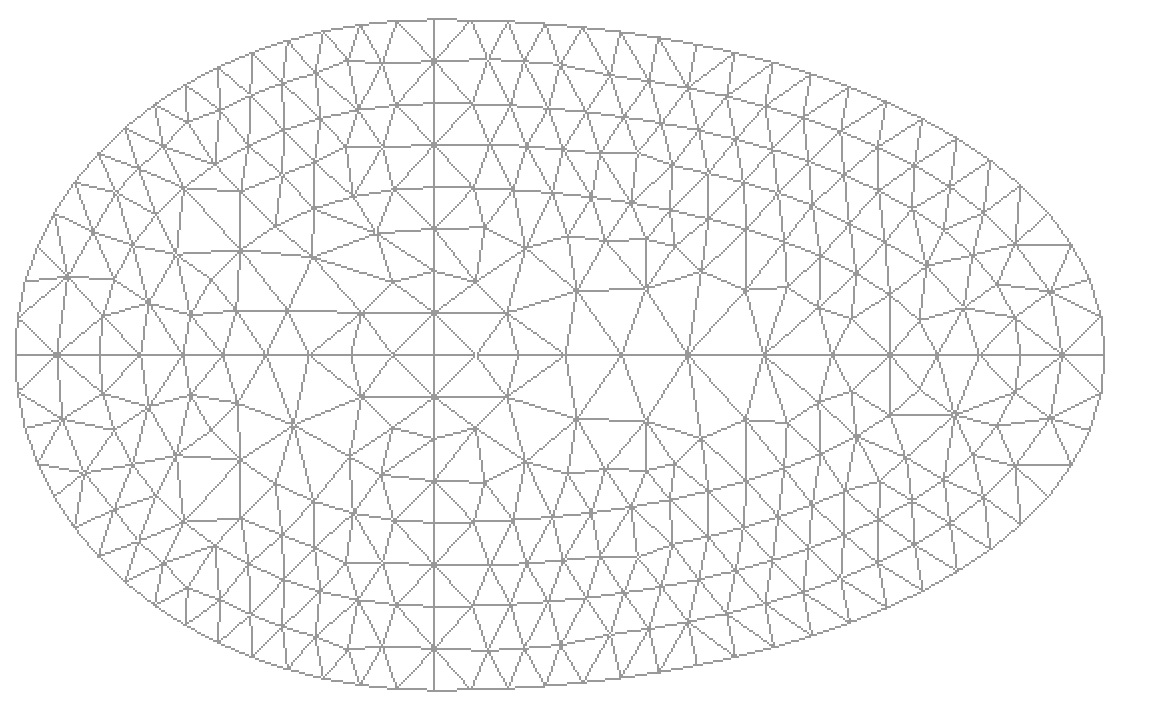}}\hspace{5mm}
	\subfigure[]{\includegraphics[width=0.35\textwidth]{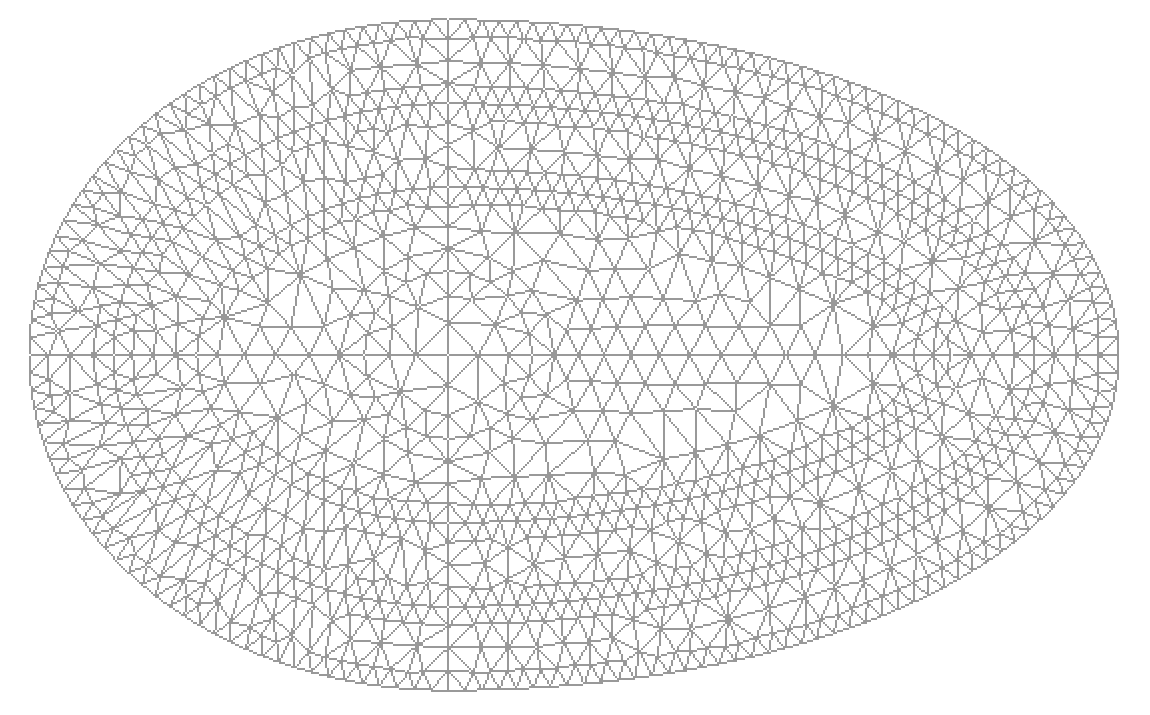}}
	\subfigure[]{\includegraphics[width=0.35\textwidth]{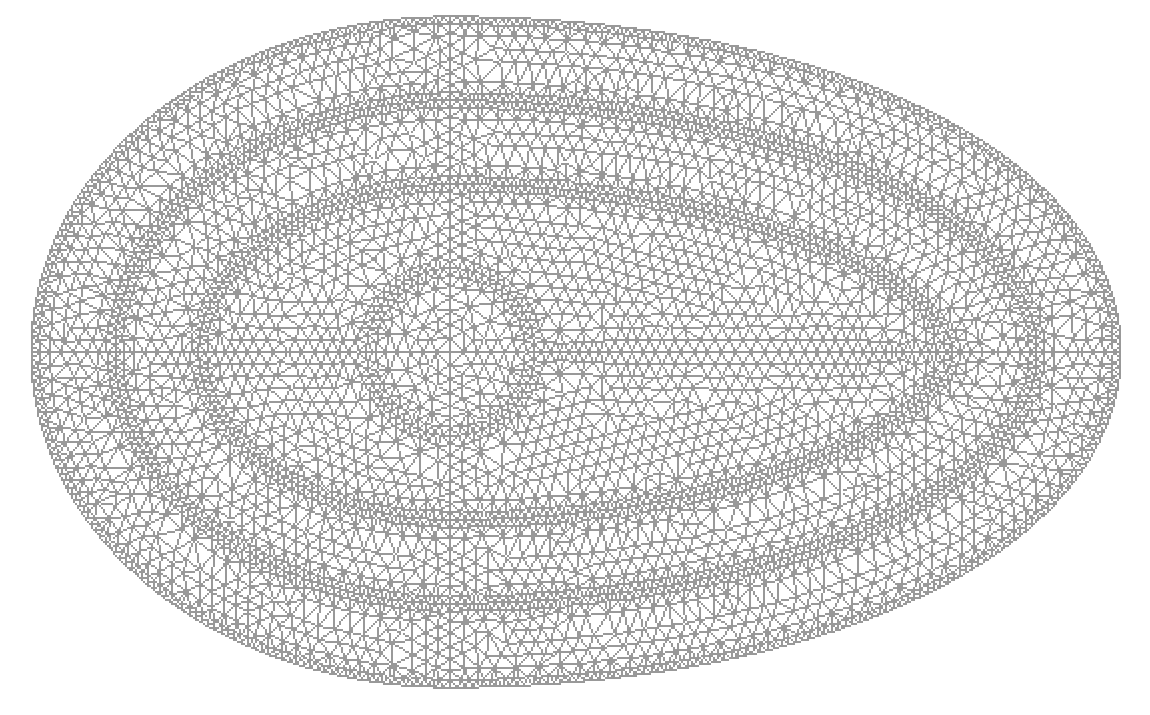}}\hspace{5mm}
	\subfigure[]{\includegraphics[width=0.35\textwidth]{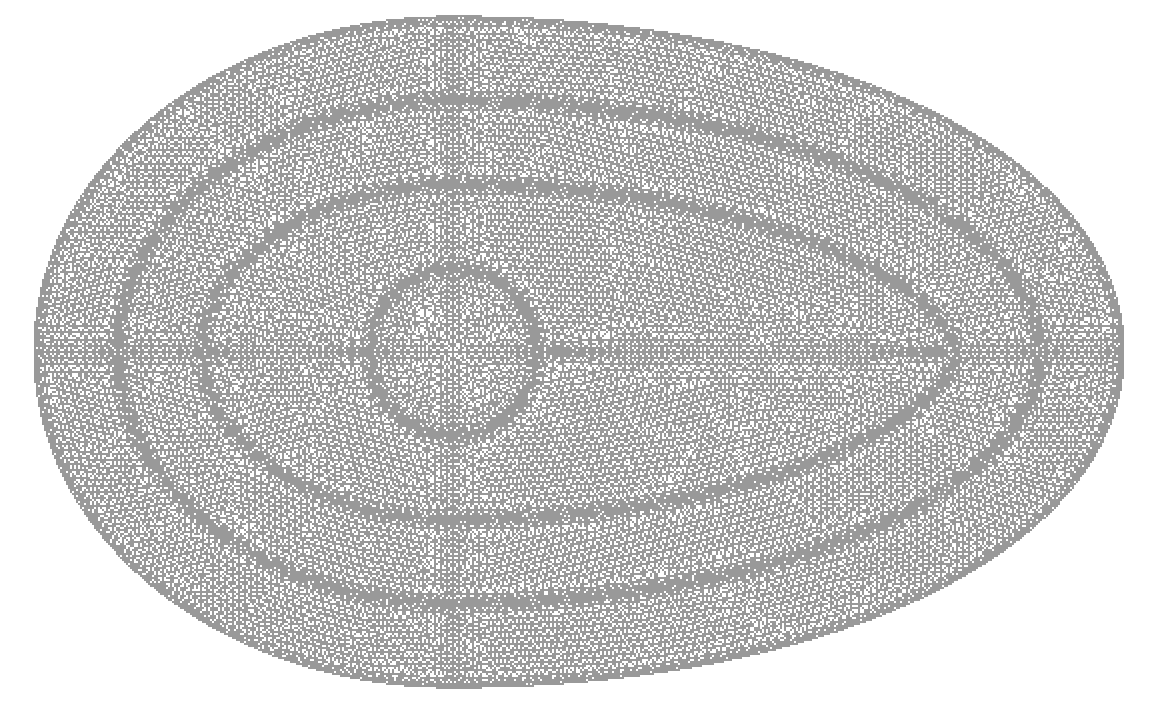}}
	\captionof{figure}{Meshes for layered egg-like figure.}\label{fig: third example meshes}
\end{figure}

\begin{figure}
	\centering
	\subfigure[Contour Fill of Temperatures]{\includegraphics[width=0.44\textwidth]{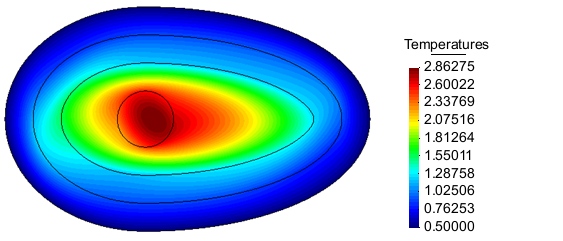}}\hspace{2mm}
	\subfigure[Contour Fill of Flux vectors on Elems]{\includegraphics[width=0.44\textwidth]{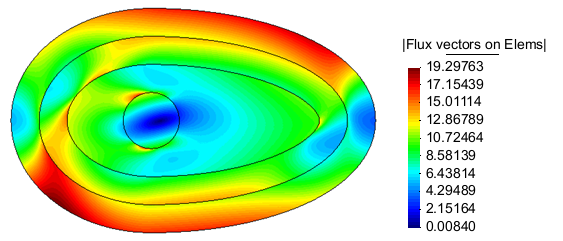}}
	\captionof{figure}{Temperature distribution and Flux magnitude fields for the finest mesh of the third example.}\label{fig: third example temperature field}
\end{figure}

\begin{figure}
	\centering
	\subfigure[]{\includegraphics[width=0.42\textwidth]{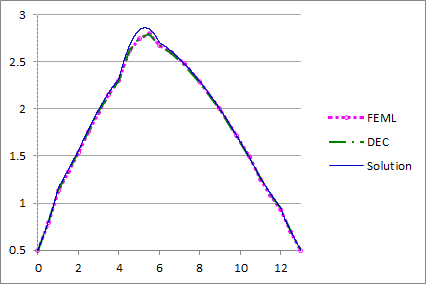}}\hspace{2mm}
	\subfigure[]{\includegraphics[width=0.42\textwidth]{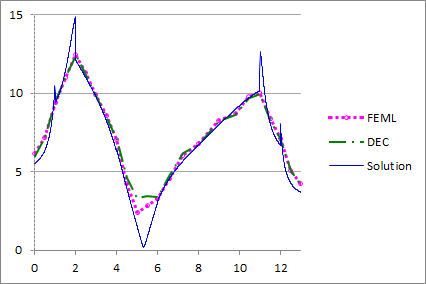}}
	\subfigure[]{\includegraphics[width=0.42\textwidth]{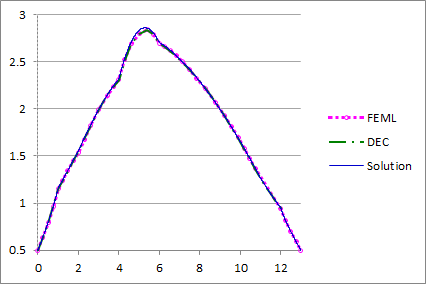}}\hspace{2mm}
	\subfigure[]{\includegraphics[width=0.42\textwidth]{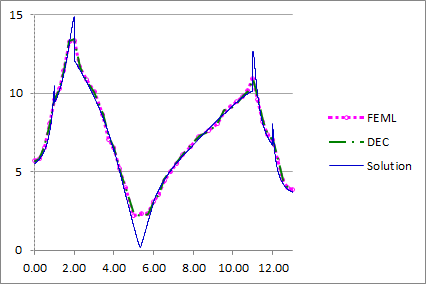}}
	\subfigure[]{\includegraphics[width=0.42\textwidth]{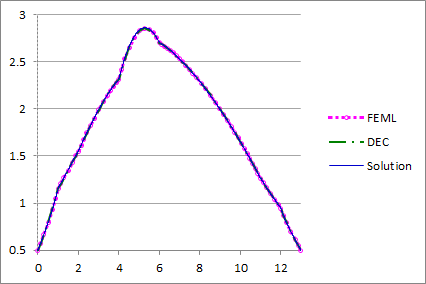}}\hspace{2mm}
	\subfigure[]{\includegraphics[width=0.42\textwidth]{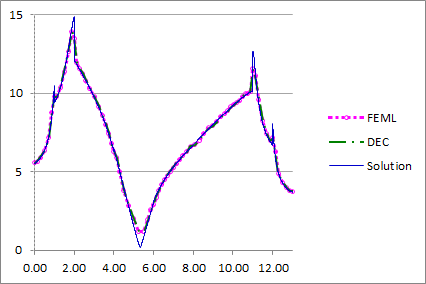}}
	\captionof{figure}{Temperature and Flux magnitude graphs of the third example along a cross-section of the domain for different meshes: 
		Mesh in Figure \ref{fig: third example meshes}(a), a-Temperature, b-Flux; 
		Mesh in Figure \ref{fig: third example meshes}(b), c-Temperature, d-Flux; 
		Mesh in Figure \ref{fig: third example meshes}(c), e-Temperature, f-Flux; 
	}\label{fig: third example diametral graphs}
\end{figure}

Let us solve the Poisson equation in a circle of radius on the following domain (see Figure \ref{fig: third example geometry}) with various material properties.
The geometry of the domain is defined by segments of ellipses passing through the given points which have centers at the origin $(0,0)$. 
\begin{center}
	\begin{tabular}{|c|c|c|c|c|c|}\hline
		Point &$x$ &$y$ &Point &$x$ &$y$ \\\hline
		a &-5 &0 &A &0 &-4 \\
		b &-4 &0 &B &0 &-3 \\
		c &-3 &0 &C &0 &-2 \\
		d &-1 &0 &D &0 &-1 \\
		e &1 &0 &E &0 &1 \\
		f &6 &0 &F &0 &2 \\
		g &7 &0 &G &0 &3 \\
		h &8 &0 &H &0 &4 \\\hline
	\end{tabular}
\end{center}
\begin{itemize}
	\item The Dirichlet boundary condition is $u=10$ and material properties (anisotropic heat diffusion constants, material angles and source terms) are given according to Figure \ref{fig: third example dirichlet condition} and the table below.
	\begin{figure}[h]
		\centering
		\includegraphics[width=0.38\textwidth]{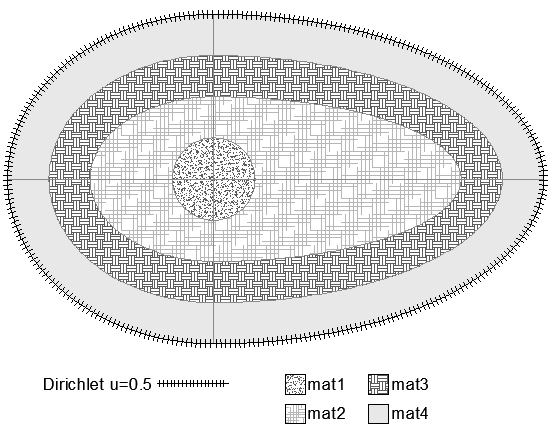}%
		\captionof{figure}{Dirichlet condition.}\label{fig: third example dirichlet condition}
	\end{figure}
	
	\begin{center} 
		\begin{tabular}{|l|c|c|c|c|}\hline
			&    $K_x$	&$K_y$&	angle&	$q$\\\hline
			Domain mat1&	5	&25	&30	&15\\\hline
			Domain mat2&	25	&5	&0	&5\\\hline
			Domain mat3&	50	&12&	45&	5\\\hline
			Domain mat4&	10	&35	&0	&5\\\hline
		\end{tabular}
	\end{center}
	
\end{itemize}
The meshes used in this example are shown in Figure \ref{fig: third example meshes}.
The numerical results for the maximum temperature value ($u(0,0)=10.2$) are exemplified in Table \ref{table: numerical results huevo} where 
a comparison with the Finite Element Method with linear interpolation functions (FEML) is also shown. 
\begin{center}
	\begin{tabular}{|c|c|c|c|c|c|c|}\hline
		\multirow{2}{*}{\rm Mesh}&\multirow{2}{*}{	\# {\rm nodes}}&\multirow{2}{*}{\# {\rm elements}}	&	\multicolumn{2}{c|}{\mbox{Max. Temp. Value}}&	\multicolumn{2}{c|}{\mbox{Max. Flux Magnitude}} \\
		&&&{\rm DEC} 	&{\rm FEML}&{\rm DEC} 	&{\rm FEML}	 \\\hline
		Figure \ref{fig: third example meshes}(a) &
		342 &
		616 &
		2.79221 &
		2.79854 &
		18.41066 &
		18.40573 \\
		Figure \ref{fig: third example meshes}(b) &
		1,259 &
		2,384 &
		2.83929 &
		2.84727 &
		18.93838 &
		18.91532 \\
		Figure \ref{fig: third example meshes}(c) &
		4,467 &
		8,668 &
		2.85608 &
		2.85717 &
		19.13297 &
		19.13193 \\
		Figure \ref{fig: third example meshes}(d)  &
		14,250  &
		28,506 &
		2.85994 &
		2.86056 &
		19.20982 &
		19.20909 \\
		&
		20,493 &
		40,316 &
		2.86120 &
		2.86177 &
		19.23120 &
		19.23457 \\
		&
		60,380 &
		119,418 &
		2.86219 &
		2.86231 &
		19.26655 &
		19.26628 \\
		&
		142,702 &
		283,162 &
		2.86249 &
		2.86256 &
		19.28045 &
		19.28028 \\
		&
		291,363 &
		579,360 &
		2.86263 &
		2.86267 &
		19.28727 &
		19.28755 \\
		&
		495,607 &
		986,724 &
		2.86275 &
		2.86269 &
		19.29057 &
		19.29081 \\
		&
		1,064,447 &
		2,122,160 &
		2.86272 &
		2.86273 &
		19.29385 &
		19.29389 \\
		&
		2,106,077 &
		4,202,536 &
		2.86274 &
		2.86274 &
		19.29618 &
		19.29615 \\
		&
		4,031,557 &
		8,049,644 &
		2.86275 &
		2.86275 &
		19.29763 &
		19.29765 \\
		\hline
	\end{tabular}    
	\captionof{table}{Maximum temperature and Flux magnitude values in the numerical simulations of the third example.}\label{table: numerical results huevo} 
\end{center}
The temperature distribution and Flux magnitude fields for the finest mesh are shown in Figure \ref{fig: third example temperature field}.
Figure \ref{fig: third example diametral graphs}
shows the graphs of the temperature and flux magnitude values along a diameter of the circle for different meshes of Figure 
\ref{fig: third example meshes}

\newpage

{\bf Remark}. As can be seen from the previous examples, DEC behaves well
on coarse meshes. 
As expected, the results of DEC and FEML are similar for fine meshes. 
We would also like to point out the the computational costs of DEC and FEML are very similar.

\section{Conclusions}\label{sec: conclusions}

DEC is a relatively recent discretization scheme for PDE's which takes into account the geometric and analytic 
features of the operators and the domains involved. 
The main contributions of this paper are the following:
\begin{enumerate}
 \item We have made explicit the local formulation of DEC, i.e. on each triangle of the mesh. As is customary, the local pieces can be assembled, which facilitates the implementation of DEC by the interested reader. 
 Furthermore, the profiles of the assembled DEC matrices are equal to those of assembled FEML matrices.

\item Guided by the local formulation,
we have deduced a natural way to approximate the flux/gradient vector of a discretized function as well as the anisotropic flux vector.
We have shown that the formulas defining the flux in DEC and FEML coincide. 

\item We have deduced how the anisotropy tensor acts on primal 1-forms.

\item We have deduced the local DEC formulation of the 2D anisotropic Poisson equation, and
have proved that the DEC and FEML diffusion terms are identical, while the source terms are not 
-- due to the different area-weight allocation for the nodes.

\item Local DEC allows a simple treatment of heterogeneous material properties assigned to subdomains (element by element), 
which eliminates the need of dealing with it through ad hoc modifications of the global discrete Hodge star operator matrix. 

\end{enumerate}

\vspace{.1in}

On the other hand we would like to point the following features:
\begin{itemize}
\item 
 The area weights assigned to the nodes of the mesh when solving the 2D anisotropic Poisson equation can even be negative (when a triangle has an inner angle greater that $120^\circ$), in stark contrast to the FEML formulation.

\item The computational cost of DEC is similar to that of FEML. While the numerical results of DEC and FEML on fine meshes are virtually identical, the DEC solutions are better than those of FEML on coarse meshes. Furthermore, DEC  solutions display numerical convergence.

\end{itemize}

\vspace{.1in}

Our future work will include the DEC discretization of convective terms and DEC on 2-dimensional simplicial surfaces in 3D. 
Preliminary results on both problems are promising and competitive with FEML. 

\bigskip

{\small
\renewcommand{\baselinestretch}{0.5}
 }


\begin{thebibliography}{10}

\bibitem{Botello} S. Botello, M.A. Moreles, E. O\~nate: ``Modulo de Aplicaciones del M\'etodo de los Elementos Finitos para resolver la ecuaci\'on de Poisson: MEFIPOISS.''  
Aula CIMNE-CIMAT, Septiembre 2010, ISBN 978-84-96736-95-5.

\bibitem{Cartan} Cartan, E.: "Sur certaines expressions diff\'erentielles et le probl\`eme de Pfaff". Annales Scientifiques de l'\'Ecole Normale Sup\'erieure. S\'erie 3. Paris: Gauthier-Villars. 16: 239?332 (1899) 

\bibitem{Craneetal} Crane, Keenan, et al. "Digital geometry processing with discrete exterior calculus." ACM SIGGRAPH 2013 Courses. ACM, 2013.

\bibitem{Dassiosetal} Dassios, Ioannis, et al. "A mathematical model for plasticity and damage: A discrete calculus formulation." Journal of Computational and Applied Mathematics 312 (2017): 27-38.

\bibitem{Esqueda1} Esqueda, H.; Herrera, R.; Botello, S; Moreles, M. A.:
"A geometric description of Discrete Exterior Calculus for general triangulations".
Rev. int. m\'etodos num\'er. c\'alc. dise\~no ing. (Online first). 

\url{https://www.scipedia.com/public/Herrera_et_al_2018b}

\bibitem{Griebel_R_S} Griebel, Michael, Christian Rieger, and Alexander Schier. "Upwind Schemes for Scalar Advection-Dominated Problems in the Discrete Exterior Calculus." Transport Processes at Fluidic Interfaces. Birkh\"auser, Cham, 2017. 145-175.

\bibitem{HiraniThesis} Hirani, Anil Nirmal. "Discrete exterior calculus". Diss. California Institute of Technology, 2003.

\bibitem{Hirani_K_N}Hirani, Anil N., Kalyana B. Nakshatrala, and Jehanzeb H. Chaudhry. "Numerical method for Darcy flow derived using Discrete Exterior Calculus." International Journal for Computational Methods in Engineering Science and Mechanics 16.3 (2015): 151-169.

\bibitem{Mohamedetal} Mohamed, Mamdouh S., Anil N. Hirani, and Ravi Samtaney. "Discrete exterior calculus discretization of incompressible Navier-Stokes equations over surface simplicial meshes." Journal of Computational Physics 312 (2016): 175-191.

\bibitem{Onate} E. O\~nate: 
``4 - 2D Solids. Linear Triangular and Rectangular Elements," in Structural Analysis with the Finite Element Method. 
Linear Statics, Volume 1: Basis and Solids, CIMNE-Springer, Barcelona, 2009. Pages 117-157, ISBN 978-1-4020-8733-2 

\bibitem{Zienkiewicz1} 
O. C. Zienkiewicz, R. L. Taylor and J. Z. Zhu: ``3 - Generalization of the finite element concepts. Galerkin-weighted residual and variational approaches,'' In The Finite Element Method Set (Sixth Edition), Butterworth-Heinemann, Oxford, 2005, Pages 54-102, ISBN 9780750664318, 


\end{thebibliography}
\end{document}